\documentclass[leqno,11pt]{article}

\usepackage{amsmath}
\usepackage{epsfig}
\usepackage[latin1]{inputenc}
\usepackage{amssymb}
\usepackage[all]{xy}
\usepackage{color}
\usepackage[normalem]{ulem}

\newcommand{\qed}{\penalty 500\hfill$\square$\par\medskip}

\newcommand{\R}{\mathrm{I\!R\!}}

\newcommand{\dd}{\mathrm{d}}

\newtheorem{teor}{Theorem}
\newtheorem{propo}{Proposition}
\newtheorem{defi}{Definition}
\newtheorem{lema}{Lemma}

\newcommand{\n}{\noindent}

\begin{document}

\title{\Large{\bf Eigenvalue for a problem involving the fractional $(p,q)$-Laplacian operator and nonlinearity with a singular and a supercritical Sobolev growth}}

\author{\sc
A. L. A. de Araujo
\thanks{Departamento de Matem\'{a}tica,
	Universidade Federal de Viçosa, 36570-900 - Vi\c{c}osa - MG, Brazil. E-mail: anderson.araujo@ufv.br}
\,\,\,
Aldo H. S. Medeiros
\thanks{Departamento de Matem\'{a}tica,
	Universidade Federal de Viçosa, 36570-900 - Vi\c{c}osa - MG, Brazil. E-mail: aldo.medeiros@ufv.br}
}

\date{}

\maketitle

\setcounter{equation}{0}
\begin{abstract}
\end{abstract}

In this paper, we are interested in studying the multiplicity, uniqueness, and nonexistence of solutions for a class of singular elliptic eigenvalue problem for the Dirichlet fractional $(p,q)$-Laplacian. The nonlinearity considered involves supercritical Sobolev growth. Our approach is variational togheter with the sub- and supesolution methods, and in this way we can address a wide range of problems not yet contained in the literature. Even when $W^{s_1,p}_0(\Omega) \hookrightarrow L^{\infty}\left(\Omega\right)$ failing, we establish $\|u\|_{L^{\infty}\left(\Omega\right)} \leq C[u]_{s_1,p}$ (for some $C>0$ ), when $u$ is a solution. 

\  

\noindent {\bf $2020$ Mathematics Subject Classification:} 35J75 35R11	35J67 35A15 

\ 


\

\date{}

\maketitle

\section{Introduction and main results}
\mbox{}

Let $\Omega \subset \mathbb{R}^{N}$ be a bounded domain. In this paper, we study the following singular eigenvalue problem for the Dirichlet fractional $(p,q)$-Laplacian
\begin{equation} \label{P0}\tag{$P_\lambda$}
\left\{
\begin{array}{llll}
(-\Delta_p)^{s_1}u + (-\Delta_q)^{s_2}u = \lambda \left[u(x)^{-\eta} + f(x,u)\right] & {\rm in} \ \ \Omega,\\
u= 0 & {\rm in} \ \R^N\setminus\Omega,\\
u> 0 & {\rm in} \ \Omega
\end{array}
\right.
\end{equation}
with $\lambda>0$, $0<s_2<s_1<1$, $0<\eta <1$ and $1<q<p$.

The fractional $p$-laplacian operator $(-\Delta_p)^{s}$ is defined as

\[
(-\Delta_p)^{s}u(x) = C(N,s,p)\lim_{\varepsilon\searrow 0}\int\limits_{\R^{N}\setminus B_\varepsilon(x)}\frac{\vert u(x) - u(y) \vert^{p-2} (u(x)-u(y))}{\vert x-y\vert^{N+sp}}\; dy\, ,
\]
\n for all $x \in \R^{n}$, where $C(N,s,p)$ is a normalization factor. The fractional p-Laplacian
is a nonlocal version of the $p$-Laplacian and is an extension of the fractional Laplacian $(p = 2)$. 

In $\left(P_\lambda\right)$, we have the sum of two such operators. So, in problem $\left(P_\lambda\right)$, the differential operator is nonhomogeneous, and this is a source of difficulties in the study of $\left(P_\lambda\right)$. Boundary value problems, driven by a combination of two or more operators of different natures, arise in many mathematical models of physical processes. One of the first such models was introduced by Cahn-Hilliard \cite{Can} describing the process of separation of binary alloys. Other applications can be found in Bahrouni-Radulescu-Repovs \cite{Radu} (on transonic flow problems).

In the reaction of $\left(P_\lambda\right), \lambda>0$ is a parameter, $u \mapsto u^{-\eta}$ with $0<\eta<1$ is a singular term and $f(z, x)$ is a Carathéodory perturbation (that is, for all $x \in \mathbb{R}, z \mapsto f(z, x)$ is measurable on $\Omega$ and for a.a $z \in \Omega, x \mapsto f(z, x)$ is continuous). Unlike many authors, we will not assume that for a.a $z \in \Omega, f(z, \cdot)$ is $(p-1)$-superlinear near $+\infty$. However, this superlinearity of the perturbation $f(z, \cdot)$ is not formulated using the very common in the literature Ambrosetti-Rabinowitz condition (the AR-condition, for short), see Ref. \cite{Sing}. The main goal of the paper is to explore the existence of a positive solution to $\left(P_\lambda\right)$. Using variational tools from the critical point theory together with truncations and comparison techniques, we show that $\left(P_\lambda\right)$ has a positive solution.

Throughout this paper, to simplify notation, we omit the constant $C(N,s,p)$. From now on, given a subset $\Omega$ of $R^N$ we set $\Omega^c = R^N\backslash \Omega$ and $\Omega^2 = \Omega \times \Omega$.  The fractional Sobolev spaces $W^{s,p}(\Omega)$ is defined to be the set of functions $u \in L^p(\Omega)$ such that
\[
\left[ u \right]_{s,p} = \left(\int\limits_{\R^{N}}\int\limits_{\R^{N}}\frac{\vert u(x) - u(y) \vert^{p}}{\vert x-y\vert^{N+sp}}\; dxdy\right)^{\frac{1}{p}} < \infty.
\]
and we defined the space $W^{s,p}_0(\Omega)$ by
\[
W^{s,p}_0(\Omega) = \bigg\{ u \in W^{s,p}(\Omega); \ \ u = 0 \ \ \text{in} \ \ \Omega^c\bigg\}.
\]

In \cite{Brasco} the authors showed that,
\[
W^{s_1,p}_0(\Omega) \hookrightarrow W^{s_2,q}_0(\Omega)
\]
thus the ideal space for to study the problem \eqref{P0} is $W^{s_1,p}_0(\Omega)$.

The main spaces that will be used in the analysis of problem \eqref{P0} are the Sobolev space $W_0^{s_1,p}(\Omega)$ and the Banach space 
\[
C^0_{s_1}(\overline{\Omega}) =\big\{ u \in C^0(\overline{\Omega}); \frac{u}{d_{\Omega}^{s_1}} \ \text{has a continuos extension to }  \overline{\Omega}\big\}.
\]
where $d_{\Omega}$ is the distance function, $d_{\Omega} =  \text{dist}(x,\partial\Omega)$.

On account of the Poincaré inequality, we have that $\left[.\right]_{s,p}$ is a
norm of the Sobolev space $W^{{s_1},p}_0(\Omega)$. Moreover,  in \cite{Brasco} the authors show that
\[
\left[ u \right]_{s_2,p} \leq \frac{C}{s_2(s_1-s_2)} \left[ u \right]_{s_1,p}, \ \ \text{for all} \ \ u \in W_0^{s_1,p}(\Omega),
\]
for  $0<s_2<s_1<1$ and $1<p<q<\infty$, in other words, we have $W_{0}^{s_1,p}(\Omega) \hookrightarrow W_0^{s_2,q}(\Omega)$.

The Banach space $C^0_{s_1}(\overline{\Omega})$  is ordered with positive (order) cone
\[
(C^0_{s_1}(\overline{\Omega}))_+ =\bigg\{f \in C^0_{s_1}(\overline{\Omega}) ; \ \ f \geq 0 \ \ \text{in} \ \ \Omega\bigg\}
\]
which is nonempty and has topological interior 
\[
\text{int} \left(C^0_{s_1}(\overline{\Omega})_+\right) = \bigg\{ v \in C^0_{s_1}(\overline{\Omega}); \ \ v > 0 \ \ \text{in} \ \ \Omega \ \ \text{and} \ \ \inf\frac{v}{d_{\Omega}^{s_1}} > 0 \bigg\}.
\]

Given $u,v \in W_0^{s_1,p}(\Omega)$ with $u \leq v$ we denote
\begin{align*}
[u,v] &= \{h \in W_0^{s_1,p}(\Omega); \ u(x) \leq h(x) \leq v(x) \ \ \text{for a. a.} \ \Omega\}\\
[u) &= \{h \in W_0^{s_1,p}(\Omega); \ \ u(x) \leq h(x) \ \text{for a. a.} \ \Omega\}.
\end{align*}
\section{The hypotheses}

The hypotheses on the perturbation $f(x,t)$ are following:
\begin{enumerate}
\item[$\textbf{H}$:] $f: \Omega \times \mathbb{R} \to \mathbb{R}$ is a Carathéodory function such that $f(x,0) = 0$ for a. a. $x \in \Omega$ and for each $t >0$ fixed $f(\cdot,t), \frac{1}{f(\cdot,t)} \in L^{\infty}(\Omega)$, moreover
	\begin{enumerate}
		\item[(i)] $\displaystyle\lim_{n\to \infty} \displaystyle\frac{F(x,t)}{t^p} = \infty$ uniformly for a. a. $x \in \Omega$, where $F(x,t) = \displaystyle\int_{0}^{t} f(x,s)\dd s$;
		\item[(ii)] if $e(x,t) = \left[1- \displaystyle\frac{p}{1-\eta}\right]t^{1-\eta} + f(x,t).t -pF(x,t)$, then there exists $\beta \in (L^{1}(\Omega))_+$ such that
		\[
		e(x,t) \leq e(x,s) + \beta(x) \ \ \text{for a.a.} \ \ x\in \Omega \ \ \text{all} \ \ 0\leq t\leq s.
		\]
		\item[(iii)] there exist $\delta>0$ and $\tau \in (1,q)$ and $c_0 >0$ such that,
		\[
		c_0t^{\tau-1} \leq f(x,t) \ \ \text{for a.a.} \ \ x \in \Omega \ \ \text{all} \ \ t \in[0,\delta]
		\]
		and for $s>0$, we have 
		\[
		0<m_s\leq f(x,t) \ \ \text{for a.a.} \ \ x \in \Omega \ \ \text{all} \ \ t\geq s.
		\]
		\item[(iv)] for every $\rho>0$, there exists $\widehat{E}_{\rho} >0$ such that for a.a. $x\in \Omega$, the function
		\[
		t\mapsto f(x,t) + \widehat{E}_{\rho}t^{p-1}
		\]
		is nondecreasing on $[0,\rho]$.
		\item[(v)]  We assume that there exists a number $\theta>0$ such that
		$$
		\limsup _{t \rightarrow \infty} \frac{f(x,t)}{t^{p^*_{s_1}-1+\theta}}<+\infty \mbox{ uniformly in } x.
		$$
		\item[(vi)] At last, we assume that there exists a sequence $\left(M_k\right)$ with $M_k \rightarrow \infty$ and such that, for each $r \in (p, p^*_{s_1})$,
		$$
		t \in\left[0, M_k\right] \Longrightarrow \frac{f(x,t)}{t^{r-1}} \leq \frac{f\left(x, M_k\right)}{\left(M_k\right)^{r-1}} \mbox{ uniformly in } x.
		$$
\end{enumerate}
\end{enumerate}
The classical AR-condition restricts $f(x,.)$ to have at least $(\mu-1)$-polynomial growth near $\infty$. In contrast, the quasimonotonicity condition that we use in this work (see hypothesis $\textbf{H} \ (ii)$), does not impose such a restriction on the growth of $f(x,.)$ and permits also the consideration of superlinear nonlinearities with slower growth near $\infty$ (see the examples below). Besides, hypothesis ($\textbf{H} \ (ii)$) is a slight extension of a condition used by Li-Yang \cite[condition $(f_4)$]{Li-Yang}.

There are convenient ways to verify ($\textbf{H} \ (ii)$). So, the hypothesis ($\textbf{H} \ (ii)$) holds, if we can find $M>0$ such that for a.a $x \in \Omega$ 
\begin{enumerate}
\item[$\bullet$] $t \mapsto \displaystyle\frac{t^{-\eta} + f(x,t)}{t^{p-1}} \ \ \text{is nondecreasing on} \ \ [M,\infty)$.
\item[$\bullet$] or $t \mapsto e(x,t) \ \ \text{is nondecreasing on} \ \ [M,\infty)$.
\end{enumerate}
	
Hypothesis ($\textbf{H} \ (iii)$) implies the presence of a concave term near zero, while hypothesis ($\textbf{H} \ (iv)$) is a one sided local Hölder condition. It is satisfied, if for a.a $x \in \Omega$, $f(x,.)$ is differentiable and for every $\rho>0$ we can find $\widehat{c}_\rho$ such that
\[
-\widehat{c}_{\rho}t^{p-1} \leq f'_t(x,t)t \ \ \text{for a.a} \ \ x \in \Omega, \ \ \text{all} \ \  0\leq t \leq \rho.
\]

Below we list two examples of functions that satisfy the conditions $(\textbf{H})$
\begin{enumerate}
	\item[$\bullet$] The function $f_1(x,t) = \left\{
	\begin{array}{llll}
		t^{\tau -1}& {\rm if} \ \ 0\leq t \leq 1,\\
		t^{p^*_{s_1}-1+\theta} & {\rm if} \ t>1,
	\end{array}
	\right.$ with $1<\tau<q<p<\theta<p_{s_1}^*$ satisfies the hipotheses $(\textbf{H}$) and also the AR-condition.
	\item[$\bullet$] The function $f_2(x,t) = \left\{
	\begin{array}{llll}
		t^{\tau -1}& {\rm if} \ \ 0\leq t \leq 1,\\
		t^{p^*_{s_1}-1+\theta}\ln{t} + t^{s-1} & {\rm if} \ t>1,
	\end{array}
	\right.$ with $1<\tau<q<p, \ \ 1 < s < p$ satisfies the hipotheses $(\textbf{H}$) but does not satisfy the AR-condition.
\end{enumerate}

\section{Preliminary}

For any $r>1$ consider the function $J_r:\mathbb{R} \to \mathbb{R}$ given by $J_r(t) = \vert t\vert^{r-2}.t$. Thus, using the arguments of \cite{Simon}, there exists $c_r>0$ and $\tilde{c}_r>0$ such that
\begin{align}
\left\langle J_r(z) - J_r(w), z-w \right\rangle &\geq \left\{
\begin{array}{ll}
c_r \vert z-w\vert^{r}, & {\rm if} \ \ r\geq 2,\\
c_r\displaystyle\frac{\vert z-w \vert^2}{\left(\vert z \vert + \vert w \vert\right)^{2-r}},  & {\rm if} \ \ r\leq 2.
\end{array}
\right. \label{eq1}\\
\vert J_r(t_1) - J_r(t_2) \vert \leq & \left\{
\begin{array}{ll}
\tilde{c_r} \vert t_1 - t_2\vert^{r-1}, & {\rm if} \ \ r\leq 2,\\
\tilde{c_r}\vert t_1-t_2\vert^2 . \left(\vert t_1 \vert + \vert t_2 \vert\right)^{r-2},  & {\rm if} \ \ r\geq 2.
\end{array}
\right.\label{eq2}
\end{align}

\begin{lema}\label{Lema1}
Let $u, v \in W_0^{s,r}(\Omega)$ and denote $w = u-v$. Then, 
\begin{equation*}
\int_{\mathbb{R}^{2N}} \frac{\big(J_r(u(x) - u(y)) - J_r(v(x) - v(y))\big)\big(w(x) - w(y)\big)}{\vert x-y \vert^{N+sr}} \dd x \dd y \geq \left\{
\begin{array}{ll}
c_r \left[ u - v \right]_{s,r}^{r}, & {\rm if} \ \ r\geq 2,\\
c_r\displaystyle\frac{\left[ u-v \right]_{s,r}^2}{\left(\left[ u \right]_{s,r} + \left[ v \right]_{s,r}\right)^{2-r}},  & {\rm if} \ \ r\leq 2.
\end{array}
\right. 
\end{equation*}
\end{lema}
\n {\bf Proof.}
The case $r \geq 2$, the result is an immediate application of the above inequality.

\textbf{Case $r\leq 2$}. Note that, using the Holder inequality we have
\begin{align*}
&\int_{\mathbb{R}^{2N}} \frac{\vert u(x) - u(y) \vert^{r}}{\vert x-y \vert^{N+sr}} \dd x \dd y = \int_{\mathbb{R}^{2N}} \frac{\vert u(x) - u(y) \vert^{r}}{\vert x-y \vert^{N+sr}}.\frac{\big(\vert u(x) - u(y) \vert + \vert v(x) - v(y) \vert\big)^{\frac{r(2-r)}{2}}} {\big(\vert u(x) - u(y) \vert + \vert v(x) - v(y) \vert\big)^{\frac{r(2-r)}{2}}}\dd x \dd y\\
&= \int_{\mathbb{R}^{2N}} \left[\frac{\vert u(x) - u(y) \vert}{\left(\vert u(x) - u(y) \vert + \vert v(x) - v(y) \vert\right)^{\frac{(2-r)}{2}}\vert x-y \vert^{\frac{N+sr}{2}}}\right]^{r}\frac{\big(\vert u(x) - u(y) \vert + \vert v(x) - v(y) \vert\big)^{\frac{r(2-r)}{2}}}{\vert x-y\vert^{\frac{2-r}{2}}}\dd x \dd y\\
&\leq \left(\int_{\mathbb{R}^{2N}} \frac{\vert u(x) - u(y) \vert^2}{\big(\vert u(x) - u(y) \vert + \vert v(x) - v(y) \vert\big)^{2-r}\vert x-y \vert^{N+sr}} \dd x \dd y\right)^{\frac{r}{2}} \left(\left[ u \right]_{s,r} + \left[ v \right]_{s,r} \right)^{\frac{r(2-r)}{2}}
\end{align*}

Thus, using the inequality \eqref{eq1} we have
\begin{align*}
\left(\frac{\left[ u-v \right]_{s,r}^r}{\left(\left[ u \right]_{s,r} + \left[ v \right]_{s,r}\right)^{\frac{r(2-r)}{2}}}\right)^{\frac{2}{r}} &\leq \int_{\mathbb{R^{2N}}} \frac{\vert u(x) - u(y) \vert^2}{\left(\vert u(x) - u(y) \vert + \vert v(x) - v(y) \vert\right)^{2-r}\vert x-y \vert^{N+sr}} \dd x \dd y \\
&\leq \frac{1}{c_r} \int_{\mathbb{R}^{2N}} \frac{\big(J_r(u(x) - u(y)) - J_r(v(x) - v(y))\big)\big((u-v)(x) - (u-v)(y)\big)}{\vert x-y \vert^{N+sr}} \dd x \dd y.
\end{align*}
\qed 

For every $1 < r < \infty$, denote by	$A_{s,r}:W_0^{s,r}(\Omega) \to \left(W_0^{s,r}(\Omega)\right)^{*}$ the nonlinear map defined by
\[
\langle A_{s,r}(u),\varphi \rangle = \int_{\mathbb{R}^{2N}} \frac{J_r(u(x)-u(y))(\varphi(x) - \varphi(y))}{\vert x-y \vert^{N+sr}}  \dd x \dd y, \ \ \text{for all} \ \ u, \varphi \in W_{0}^{s,r}(\Omega).
\]

An immediate consequence of Lemma \ref{Lema1} is the following proposition 	
\begin{propo}
The map $A_{s,r}:W_0^{s,r}(\Omega) \to \left(W_0^{s,r}(\Omega)\right)^{*}$  maps bounded sets to bounded sets, is continuous, strictly monotone and satisfies,
\[
u_n \rightharpoonup u \ \ \text{in} \ \ W_0^{s,r}(\Omega) \ \text{and} \ \limsup_{n\to \infty} \ \langle A_{s,r}(u_n),(u_n-u) \rangle \leq 0 \Rightarrow u_n \to u \ \ \text{in} \ \ W_0^{s,r}(\Omega).
\]
\end{propo}


\n {\bf Proof.} Indeed, using the inequality \eqref{eq2} we have
\[
\Vert A_{s,r}(u) - A_{s,r}(w) \Vert_{*} \leq \left\{
\begin{array}{ll}
\tilde{c_r} \left[ u - w \right]_{s,r}^{r-1}, & {\rm if} \ \ r\leq 2,\\
\tilde{c_r}\left[ u-w\right]_{s,r}^2 . \left(\left[ u \right]_{s,r} + \left[w \right]_{s,r}\right)^{r-2},  & {\rm if} \ \ r\geq 2.
\end{array}
\right.
\]
and thus $A_{s,r}$ maps bounded sets to bounded sets, is continuous.	

Moreover, if $p\geq 2$ then using also the Lemma \ref{Lema1} results,
\begin{align*}
\lim_{n \to \infty} c_r\left[ u_n - u \right]_{s,r}^2
&\leq \lim_{n\to \infty}\int_{\mathbb{R}^{2N}} \frac{\big(J_r(u_n(x) - u_n(y)) - J_r(u(x) - u(y))\big)\big((u_n-u)(x) - (u_n-u)(y)\big)}{\vert x-y \vert^{N+sr}} \dd x \dd y \\
&= \limsup_{n\to \infty}\bigg\langle A_{s,r}(u_n) - A_{s,r}(u), u_n-u \bigg\rangle \leq 0,
\end{align*}
and if $p\leq 2$ let's use again the Lemma \ref{Lema1} and obtain
\begin{align*}
c_r\frac{\left[ u_n - u \right]_{s,r}^2}{\left(\left[ u_n \right]_{s,r} + \left[ u \right]_{s,r}\right)^{2-r}}
&\leq \int_{\mathbb{R}^{2N}} \frac{\big(J_r(u_n(x) - u_n(y)) - J_r(u(x) - u(y))\big)\big((u_n-u)(x) - (u_n-u)(y)\big)}{\vert x-y \vert^{N+sr}} \dd x \dd y \\
&= \bigg\langle A_{s,r}(u_n) - A_{s,r}(u), u_n-u \bigg\rangle
\end{align*}
thus, if $u_n \rightharpoonup u$ in $W^{s,r}_{0}(\Omega)$ and $\displaystyle\limsup_{n\to \infty} A_{s,r}(u_n).(u_n-u) \leq 0$ then, there exists $M>0$ such that $\Vert u_n \Vert_{s,r} \leq M$ and thus
\begin{align*}
\lim_{n\to \infty} c_r\frac{\left[u_n - u\right]_{s,r}^2}{\left(M + \left[ u \right]_{s,r}\right)^{2-r}} \leq \lim_{n\to \infty}c_r\frac{\left[ u_n - u\right]_{s,r}^2}{\left(\left[ u_n \right]_{s,r} + \left[ u \right]_{s,r}\right)^{2-r}}
&\leq \limsup_{n\to \infty}\bigg\langle A_{s,r}(u_n) - A_{s,r}(u), u_n - u \bigg\rangle \leq 0.
\end{align*}
Consequently, for all $1<p<\infty$, we have $u_n \to u$ in $W_{0}^{s,r}(\Omega)$.

\qed

The following result is a natural improvement of \cite[Lemma 9]{Lind} to the Dirichlet fractional $(p,q)$- Laplacian.

\begin{propo}[Weak comparison principle]\label{PC}
Let $0<s_1<s_2<1$, $1<q<p$, $\Omega$ be bounded in $\mathbb{R}^{N}$ and $u,v \in W_0^{s_1,p}(\Omega) \cap C_{s_1}^{0}(\overline{\Omega})$. Suppose that, 
\begin{align*}
\bigg\langle A_{s_1,p}(u) + A_{s_2,q}(u), (u-v)^+ \bigg\rangle  \leq \bigg\langle A_{s_1,p}(v) + A_{s_2,q}(v),(u-v)^+ \bigg\rangle 
\end{align*}
then $u \leq v$. 
\end{propo}

\n {\bf Proof.}	The proof is a straightforward calculation, but for convenience of the reader we present a sketch of it. By considering the equations for both $p$ and $q$, and subtracting them and adjusting the terms, we obtain
\begin{align}\label{int1}
\bigg\langle A_{s_1,p}(u) + A_{s_2,q}(u), (u-v)^+ \bigg\rangle  - \bigg\langle A_{s_1,p}(v) + A_{s_2,q}(v),(u-v)^+ \bigg\rangle \leq 0.
\end{align}

Using the identity
\begin{align*}
J_{m}(b) - J_m(a) = (m-1)(b-a)\displaystyle\int_{0}^{1} \vert a + t(b-a) \vert^{m-2} \dd t
\end{align*}
for $a = v(x) - v(y)$ and $b = u(x)- u(y)$, we have
\begin{align*}
J_{m}(u(x) - u(y)) - J_{m}(v(x) - v(y)) = (m-1)\left[(u-v)(x) - (u-v)(y) \right]Q_{m}(x,y),
\end{align*}
where $Q_{m}(x,y) = \displaystyle\int_{0}^{1}\left\vert (v(x) - v(y)) + t[(u-v)(x) - (u-v)(y)]\right\vert^{m-2}\dd t$.
	
We have $Q_{m}(x,y) \geq 0$ and $Q_{m}(x,y) = 0$ only if $v(x) = v(y)$ and $u(x) = u(y)$. Rewriting the integrands in \eqref{int1} we obtain
\begin{multline*}
\int_{\mathbb{R}^{2N}} \left(\frac{(p-1)\left[(u-v)(x) - (u-v)(y) \right]Q_{p}(x,y)}{\vert x-y\vert^{N+sp}}\right)((u-v)^+(x)-(u-v)^+(y))\dd x\dd y\\
+ \int_{\mathbb{R}^{2N}} \left(\frac{(q-1)\left[(u-v)(x) - (u-v)(y) \right]Q_{q}(x,y)}{\vert x-y\vert^{N+sq}}\right)((u-v)^+(x)-(u-v)^+(y)) \dd x\dd y  \leq 0. 
\end{multline*}
We now consider 
\[
\psi = u-v = (u-v)^{+} - (u-v)^{-} , \quad \varphi = (u-v)^{+} = \psi^{+}.
\]
It follows from the last inequality that
\begin{align*}
\displaystyle\int_{\mathbb{R}^{2N}} &\left(\frac{(p-1)(\psi(x) - \psi(y) )(\psi^{+}(x)-\psi^{+}(y))Q_{p}(x,y)}{\vert x-y\vert^{N+sp}}\right)\dd x\dd y \\
&+\int_{\mathbb{R}^{2N}} \left(\frac{(q-1)(\psi(x) - \psi(y) )(\psi^{+}(x)-\psi^{+}(y))Q_{q}(x,y)}{\vert x-y\vert^{N+sq}}\right)\dd x\dd y \leq 0.
\end{align*}
Applying the inequality $(\xi - \eta)(\xi^{+} - \eta^{+}) \geq \vert \xi^{+} - \eta^{+} \vert^{2}$ we obtain
\begin{align*}
&\int_{\mathbb{R}^{2N}} \frac{(p-1)\vert \psi^{+}(x)-\psi^{+}(y)\vert^{2}Q_{p}(x,y)}{\vert x-y\vert^{N+sp}} \dd x\dd y + \int_{\mathbb{R}^{2N}}\frac{(q-1)\vert \psi^{+}(x)-\psi^{+}(y)\vert^{2}Q_{q}(x,y)}{\vert x-y\vert^{N+sq}} \dd x\dd y \leq 0.
\end{align*}

Thus, at almost every point $(x,y)$ we have $\psi^{+}(x) = \psi^{+}(y)$ or
\[
Q_{p}(x,y) = Q_{q}(x,y) = 0 .
\]

Since $Q_{p}(x,y) = Q_{q}(x,y) = 0$ also imply $\psi^{+}(x) = \psi^{+}(y)$, we conclude that
\[
(u-v)^{+}(x) = C \geq 0,  \ \ \forall x \in \mathbb{R}^{N}
\]
and since, $u,v \in W_{0}^{s_1,p}(\Omega)$, results that $C=0$ and consequently $u\leq v$.
\qed
\begin{propo}[Strong comparison principle]\label{SPC}
Let $0<s_1<s_2<1$, $1<q<p$, $\Omega$ be bounded in $\mathbb{R}^{N}$, $g \in C^{0}(\mathbb{R}) \cap BV_{loc}(\mathbb{R})$, $u,v \in W_0^{s_1,p}(\Omega) \cap C_{s_1}^{0}(\overline{\Omega})$ such that $u \neq v$ and $K>0$ satisfy, 
\begin{align*}
\left\{
\begin{array}{llll}
(-\Delta_p)^{s_1}u + (-\Delta_p)^{s_1}u + g(u)  \leq (-\Delta_p)^{s_1}v + (-\Delta_q)^{s_2}v + g(v) \leq K \ \ & {\rm weakly \  in} \ \ \Omega,\\ 
0<u\leq v \ \ &{\rm in} \ \ \Omega.
\end{array}
\right.
\end{align*}
then $u \leq v$ in $\Omega$. In particular, if $u,v \in \text{int}[(C_{s_1}^{0}(\overline{\Omega})^+)]$ then $v-u \in \text{int}[(C_{s_1}^{0}(\overline{\Omega})^+)]$.
\end{propo}
\n {\bf Proof.} Without loss of generality, we may assume that $g$ is nondecreasing
and $g(0) = 0$. In fact, by Jordan's  decomposition we can find $g_1,g_2 \in C^{0}(\mathbb{R})$ nondecreasing such that $g(t) = g_1(t) - g_2(t)$ and $g_1(0) = 0$.

Since, $u \neq v$ by continuity, we can find $x_0 \in \Omega$, $\rho,\varepsilon > 0$ such that $\overline{B_{\rho}(x_0)} \subset \Omega$ and 
\[
\sup_{\overline{B_\rho(x_0)}} u < \inf_{\overline{B_\rho(x_0)}}v - \varepsilon.
\]

Hence, for all $\eta > 1$ close enough to $1$ we have
\[
\sup_{\overline{B_\rho(x_0)}} \eta u < \inf_{\overline{B_\rho(x_0)}}v - \frac{\varepsilon}{2}.
\]
Define $w_\eta \in W_{0}^{s_1,p}(\Omega\backslash \overline{B_{\rho}(x_0)})$ by
\[
w_\eta(x) =\left\{
 \begin{array}{ll}
\eta u(x), & {\rm if} \ \ x \in \overline{B_{\rho}(x_0)}^{c},\\
v(x), & {\rm if} \ \ x \in \overline{B_{\rho}(x_0)},
\end{array}
\right.
\]
so $w_\eta \leq v(x)$ in $\overline{B_{\rho}(x_0)}$ and  by the nonlocal
superposition principle [\cite{Ian}, Proposition 2.6] we have weakly in $\Omega \backslash \overline{B_{\rho}(x_0)}$
\[
(-\Delta_p)^{s_1} w_\eta \leq \eta^{p-1} (-\Delta_p)^{s_1} u - C_{\rho}\varepsilon^{p-1} \ \ \text{and} \ \ (-\Delta_q)^{s_2} w_\eta \leq \eta^{q-1} (-\Delta_q)^{s_2} u - C_{\rho}\varepsilon^{q-1}
\]
for some $C_{\rho}>0$ and all $\eta > 1$ close enough to $1$. Further, we have weakly in $\Omega \backslash \overline{B_{\rho}(x_0)}$
\begin{align*}
(-\Delta_p)^{s_1} w_\eta &+ (-\Delta_q)^{s_2} w_\eta + g(w_\eta) \leq \eta^{p-1} (-\Delta_p)^{s_1} u + \eta^{q-1}(-\Delta_q)^{s_2} u + g(w_\eta)- C_{\rho}\varepsilon^{q-1} - C_{\rho}\varepsilon^{p-1} \\
&\leq \eta^{p-1}\bigg( (-\Delta_p)^{s_1} u + (-\Delta_q)^{s_2} u + g(u)\bigg) + \bigg(g(w_\eta) - \eta^{p-1}g(u)\bigg)-C_{\rho}\varepsilon^{q-1} - C_{\rho}\varepsilon^{p-1}\\
&\leq \bigg( (-\Delta_p)^{s_1} u + (-\Delta_q)^{s_2} u + g(u)\bigg) + \bigg(g(w_\eta) - \eta^{p-1}g(u)\bigg) + K\big(\eta^{p-1} - 1\big) - C_{\rho}\varepsilon^{q-1} - C_{\rho}\varepsilon^{p-1}\\
&\leq \bigg( (-\Delta_p)^{s_1} v + (-\Delta_q)^{s_2} v + g(v)\bigg) + \bigg(g(w_\eta) - \eta^{p-1}g(u)\bigg) + K\big(\eta^{p-1} - 1\big) - C_{\rho}\varepsilon^{q-1} - C_{\rho}\varepsilon^{p-1}.
\end{align*}

Since 
\[
\bigg(g(w_\eta) - \eta^{p-1}g(u)\bigg) + K\big(\eta^{p-1} - 1\big) \to 0
\]
uniformly in $\Omega \backslash \overline{B_{\rho}(x_0)}$ as $\eta \to 1^{+}$, we have, for all $\eta >1$ close enough to $1$,
\begin{align*}
\left\{
\begin{array}{llll}
(-\Delta_p)^{s_1}w_{\eta} + (-\Delta_p)^{s_1}w_{\eta} + g(w_\eta)  \leq (-\Delta_p)^{s_1}v + (-\Delta_q)^{s_2}v + g(v) \leq K \ \ & {\rm weakly \  in} \ \ \Omega \backslash \overline{B_{\rho}(x_0)},\\ 
0<w_{\eta}\leq v \ \ &{\rm in} \ \ \left(\Omega \backslash \overline{B_{\rho}(x_0)}\right).
\end{array}
\right.
\end{align*}

Testing with $\varphi = (w_{\eta} - v)^{+} \in W_{0}^{s_1}(\Omega)\backslash\overline{B_{\rho}(x_0)}$, recalling the monotonicity of $g$, and applying
Proposition \ref{PC} we get $v > w_\eta$ in $ \Omega)\backslash\overline{B_{\rho}(x_0)}$. So we have 
\[
v \geq \eta u \geq u.
\]
In particular, if $u, v \in \text{int}\left[(C_{s_1}^{0}(\overline{\Omega}))_+\right]$ then 
\[
\inf_{\Omega} \frac{v-u}{d^{s_1}_{\Omega}} \leq \inf_{\Omega}\frac{(\eta-1)u}{d^{s_1}_{\Omega}} > 0
\]
and so $v-u \in \text{int}\left[(C_{s_1}^{0}(\overline{\Omega}))_+\right]$. 
\qed

\section{An auxiliary problem}
Firstly, we will need to define, with the help of the real sequence defined in \textbf{H}(vii), a sequence of auxiliary equations that will be important for our purpose. More specifically, for each $k \in \mathbb{N}$, we define the auxiliary truncation functions by choosing $r \in\left(p,p^*_{s_1}\right)$ such that $p^*_{s_1}-r<\theta$ and we set
\begin{equation}\label{aprox.f}
f_k(x,t)=\left\{\begin{array}{ll}
0, & \text{if}\ t \leq 0 \\
f(x,t), & \text{if} \ 0 \leq t \leq M_k \\
\displaystyle\frac{f\left(x,M_k\right)}{\left(M_k\right)^{r-1}} t^{r-1}, & \ \ \text{if} \ \ t \geq M_k.
\end{array}\right.
\end{equation}

Notice that we define $f_k$ to be such that $r$ in its definition is independent of $k$. We see that we are really truncating our original function, making it subcritical for large arguments. Furthermore, in view of conditions \textbf{H}(vi), \textbf{H}(vii) and the choice of $\theta$, we can prove that, for $k$ big enough, $f_k$ satisfies, for a constant $C>0$,
\begin{equation}\label{est.1}
\left|f_k(x,t)\right| \leq C\left(M_k\right)^{2 \theta}|t|^{r-1}.
\end{equation}

Indeed, for all $t>0$, condition \textbf{H}(vii) and (\ref{aprox.f}) gives
$$
f_k(x,t) \leq \frac{f\left(x,M_k\right)}{\left(M_k\right)^{r-1}} t^{r-1}
$$
and, by \textbf{H}(vi), if $k$ is sufficiently large,
$$
\frac{f\left(x,M_k\right)}{\left(M_k\right)^{r-1}} \leq C\left(M_k\right)^{p^*_{s_1}-r+\theta} \leq C\left(M_k\right)^{2 \theta}.
$$
For each $k \in \mathbb{N}$, let us consider the following auxiliary problem
\begin{equation} \label{P1}\tag{$P_{k,\lambda}$}
\left\{
\begin{array}{llll}
(-\Delta_p)^{s_1}u + (-\Delta_q)^{s_2}u = \lambda \left[u(x)^{-\eta} + f_k(x,u)\right] & {\rm in} \ \ \Omega,\\
u= 0 & {\rm in} \ \R^N\setminus\Omega,\\
u> 0 & {\rm in} \ \Omega
\end{array}
\right.
\end{equation}
with $\lambda>0$, $0<s_1<s_2<1$, $0<\eta <1$ and $1<q<p$.

By the hypotheses  $(\textbf{H})$, the hypotheses on the truncation $f_k(x,t)$ are following:
\begin{enumerate}
\item[$\textbf{H}_k$:] $f_k: \Omega \times \mathbb{R} \to \mathbb{R}$ is a Carathéodory function such that $f_k(x,0) = 0$ for a. a. $x \in \Omega$ and 
\begin{enumerate}
\item[(i)] $f_k(x,t) \leq \alpha_k(x) [1 + t^{r -1}]$ for a. a. $x \in \Omega$ all $t \geq 0$ with $\alpha_k \in L^{\infty}(\Omega)$ and $p < r < p^*_{s_1} = \displaystyle\frac{NP}{N-s_1p}$;
\item[(ii)] $\displaystyle\lim_{t\to \infty} \displaystyle\frac{F_k(x,t)}{t^p} = \infty$ uniformly for a. a. $x \in \Omega$, where $F_k(x,t) = \displaystyle\int_{0}^{t} f_k(x,s)\dd s$;
\item[(iii)] if $e_k(x,t) = \left[1- \displaystyle\frac{p}{1-\eta}\right]t^{1-\eta} + f_k(x,t).t -pF_k(x,t)$, then there exists $\beta_k \in (L^{1}(\Omega))_+$ such that
\[
e_k(x,t) \leq e_k(x,s) + \beta_k(x) \ \ \text{for a.a.} \ \ x\in \Omega \ \ \text{all} \ \ 0\leq t\leq s.
\]
\item[(iv)] there exist $\delta>0$ and $\tau \in (1,q)$ and $c_0 >0$ such that,
\[
		c_0t^{\tau-1} \leq f_k(x,t) \ \ \text{for a.a.} \ \ x \in \Omega \ \ \text{all} \ \ t \in[0,\delta]
		\]
		and for all $s>0$, we have 
		\[
		0<m_{k,s}\leq f_k(x,t) \ \ \text{for a.a.} \ \ x \in \Omega \ \ \text{all} \ \ t\geq s.
		\]
		\item[(v)] for every $\rho>0$, there exists $\widehat{E}_{k,\rho} >0$ such that for a.a. $x\in \Omega$, the function
		\[
		t\mapsto f_k(x,t) + \widehat{E}_{k,\rho}t^{p-1}
		\]
		is nondecreasing on $[0,\rho]$.
	\end{enumerate}

\end{enumerate}

The hypothesis ($\textbf{H}_k \ (i)$) holds by (\ref{est.1}), ($\textbf{H}_k \ (ii)$) holds by (\ref{aprox.f}) and $p<r$. We will prove first that ($\textbf{H}_k \ (iv)$) holds. Since $\delta>0$, $\tau \in (1,q)$ and $c_0 >0$, if $\delta <M_k$, we have
\[
c_0t^{\tau-1} \leq f(x,t)=f_k(x,t) \ \ \text{for a.a.} \ \ x \in \Omega \ \ \text{all} \ \ t \in[0,\delta].
\]
For $s>0$, we have

$\bullet$ $0<s\leq t \leq M_k$,
\[
f_k(x,t)=f(x,t)\geq m_s>0,
\]
by ($\textbf{H} \ (iii)$).

$\bullet$ $0<s\leq M_k <t$,
\[
f_k(x,t)=\frac{f\left(x,M_k\right)}{\left(M_k\right)^{r-1}} t^{r-1}\geq \frac{f\left(x,M_k\right)}{\left(M_k\right)^{r-1}} M_k^{r-1} = f\left(x,M_k\right) >0.
\]

$\bullet$ $0< M_k <s\leq t$,
\[
f_k(x,t)=\frac{f\left(x,M_k\right)}{\left(M_k\right)^{r-1}} t^{r-1}\geq \frac{f\left(x,M_k\right)}{\left(M_k\right)^{r-1}} M_k^{r-1} = f\left(x,M_k\right) >0.
\]
So, for all $s>0$ we have
\[
f_k(x,t)\geq m_{k,s}>0  \ \ \text{for a.a.} \ \ x \in \Omega \ \ \text{all} \ \ t\geq s,
\]
with $m_{k,s}=\max\left\{m_s, \displaystyle\inf_{x \in \Omega}f(x,M_k)\right\}>0$.

To prove that ($\textbf{H}_k \ (iii)$) holds it is sufficiently verify that there is a constant $C_k>0$ such that $t \mapsto e_k(x,t) \ \ \text{is nondecreasing on} \ \ [C_k,\infty)$. Since for $t\geq M_k$ we have
\[
\begin{array}{rcl}
e_k(x,t) &=& \left[1- \displaystyle\frac{p}{1-\eta}\right]t^{1-\eta} + f_k(x,t).t -pF_k(x,t)\\
&=&\displaystyle \left[1- \frac{p}{1-\eta}\right]t^{1-\eta}  + \frac{f\left(x,M_k\right)}{\left(M_k\right)^{r-1}} t^{r} - p\int_0^{M_k}f(x,s)ds - \int_{M_k}^t\frac{f\left(x,M_k\right)}{\left(M_k\right)^{r-1}} s^{r-1}ds \\
&=&\displaystyle \left[1- \frac{p}{1-\eta}\right]t^{1-\eta}  + \frac{f\left(x,M_k\right)}{\left(M_k\right)^{r-1}} t^{r} - p\int_0^{M_k}f(x,s)ds - \frac{f\left(x,M_k\right)}{\left(M_k\right)^{r-1}}\frac{1}{r}[t^r-M_k^r]. 
\end{array}
\]
Hence
\[
\begin{array}{rcl}
	\displaystyle\frac{\partial}{\partial t}e_k(x,t) &=& \displaystyle\left[1-\eta- p\right]t^{-\eta} + (r-1)\frac{f\left(x,M_k\right)}{\left(M_k\right)^{r-1}} t^{r-1}. 
\end{array}
\]
Notice that $\frac{\partial}{\partial t}e_k(x,t)\geq 0$ if 
\[
\left[1-\eta- p\right]t^{-\eta} + (r-1)\frac{f\left(x,M_k\right)}{\left(M_k\right)^{r-1}} t^{r-1}\geq 0,
\]
or equivalently, if
\[
t\geq \left(-\left[1-\eta- p\right]\frac{\left(M_k\right)^{r-1}}{(r-1)f(x,M_k)}\right)^{\frac{1}{r+\eta}}.
\]
We can consider 
\[
C_k=\left(-\left[1-\eta- p\right]\frac{\left(M_k\right)^{r-1}}{(r-1)m_{k,s}}\right)^{\frac{1}{r+\eta}},
\]
where $m_{k,s}$ is as in ($\textbf{H}_k \ (iv)$). Hence, $t \mapsto e_k(x,t) \ \ \text{is nondecreasing on} \ \ [C_k,\infty)$. The proof of ($\textbf{H}_k \ (v)$) follows from \eqref{aprox.f} and ($\textbf{H} \ (iv)$).

\begin{defi}
A function $u \in W^{s_1,p}_0(\Omega)$ is a weak solution of the problem \eqref{P1} if, $u^{-\eta}\varphi \in W_0^{s_1,p}(\Omega)$ for all $\varphi \in W_0^{s_1,p}(\Omega)$ and
\[
\bigg\langle A_{s_1,p}(u) + A_{s_2,q}(u),\varphi \bigg\rangle = \int_{\Omega} \lambda \left[u^{-\eta} + f_k(x,u)\right] \varphi\dd x, \ \ \text{for all} \ \ \varphi \in W_{0}^{s_1,p}(\Omega).
\]
\end{defi}

The difficulty that we encounter in the analysis of problem \eqref{P1} is that the energy (Euler) function of the problem $I_{\lambda}: W_{0}^{s_1,p}(\Omega) \to \mathbb{R}$ defined by
\begin{equation}
I_\lambda(u) = \frac{1}{p} \left[ u \right]_{s_1,p}^p + \frac{1}{q} \left[ u \right]_{s_2,p}^{q} - \lambda \int_{\Omega} \left[\frac{1}{1-\eta}(u^{+})^{1-\eta} + F_k(x,u^+)\right] \dd x.
\end{equation}
for all $u \in W_{0}^{s_1,p}(\Omega)$, is not $C^{1}$ (due to the singular term). So, we can not use the minimax
methods of critical point theory directly on $I_\lambda(.)$. We have to find ways to bypass the singularity and to deal with $C^1$-functionals.

The hypotheses $\textbf{H} \ (i)$ and $\textbf{H} \ (iv)$ assure us that, there are $c_0 > 0$ and $c_2 > 0$ such that,
\begin{equation}\label{fcond}
f_k(x,z) \geq c_0z^{\tau - 1} - c_2z^{\theta-1}, \ \ \text{for a. a.} \ \ x \in \Omega \ \ \text{and} \ \ z \geq 0.
\end{equation}

We consider the following auxiliary Dirichilet fractional $(p,q)$-equation
\begin{equation} \label{P2}
\left\{
\begin{array}{llll}
(-\Delta_p)^{s_1}u + (-\Delta_q)^{s_2}u = \lambda \left[c_0u(x)^{\tau -1} - c_2u^{\theta -1}\right] & {\rm in} \ \ \Omega,\\
u= 0 & {\rm in} \ \R^N\setminus\Omega,\\
u> 0 & {\rm in} \ \Omega
\end{array}
\right.
\end{equation}
with $0<s_2<s_1$, $\lambda > 0$ and $1<\tau<q<p<\theta < p_s^*= \displaystyle\frac{Np}{N-sp}$.

\begin{lema}
If $\underline{u}_{\lambda}\in W_{0}^{s_1,p}(\Omega)$ be a weak solution of problem \eqref{P2}. Then $\underline{u}_{\lambda}\in L^\infty(\Omega)$.
\end{lema}	
\n {\bf Proof.}  
We denote by $h_{\lambda}(t) = \lambda c_0 t^{\tau - 1} - \lambda c_2t^{\theta-1}$. Thus, 
\begin{equation}\label{sol}
\begin{aligned}\langle A_{s_1,p}(\underline{u}_{\lambda}) + A_{s_2,q}(\underline{u}_{\lambda}),\phi \rangle &= \int_{\mathbb{R}^{2N}}\left(\frac{J_p(\underline{u}_{\lambda}(x)-\underline{u}_{\lambda}(y))}{|x-y|^{N+s_1p}} + \frac{J_q(\underline{u}_{\lambda}(x)-\underline{u}_{\lambda}(y)) }{|x-y|^{N+s_2q}}\right) (\phi(x) - \phi(y))\dd x\dd y \\
&=\int_\Omega h_{\lambda}(\underline{u}_{\lambda})\phi \dd x
\end{aligned}
\end{equation}
for any $\phi\in W_{0}^{s_1,p}(\Omega)$.

For each $k\in\mathbb{N}$, set
$$\Omega_k:=\{x\in\Omega~:~u(x)>k\}.$$

Since $\underline{u}_{\lambda}\in W_{0}^{s_1,p}(\Omega)$ and $\underline{u}_{\lambda} \geq 0$ in $\Omega$, we have that $(\underline{u}_{\lambda}-k)^+\in W_{0}^{s_1,p}(\Omega)$. Taking $\phi= (\underline{u}_{\lambda} - k)^+$ in \eqref{sol}, we obtain
\begin{equation}\label{eq 4.18}
	\left\langle A_{s_1,p}(\underline{u}_{\lambda}) + A_{s_2,q}(\underline{u}_{\lambda}),\phi \right\rangle =\int_\Omega h_{\lambda}(\underline{u}_{\lambda})(\underline{u}_{\lambda}-k)^+\dd x.
\end{equation}

Applying the  algebraic inequality $|a-b|^{p-2}(a-b)(a^+-b^+)\geq|a^+-b^+|^p$
to estimate the left-hand side of \eqref{eq 4.18}, we obtain
\begin{align}\label{eq1}
	\left(  \int_{\Omega_k}(\underline{u}_{\lambda}-k)^{p^*_s}\dd x\right)^{\frac{p}{p^*_s}} & \leq C\int_{\mathbb{R}^{2N}}\frac{|\underline{u}_{\lambda}(x)-\underline{u}_{\lambda}(y)|^p}{|x-y|^{N+sp}}\dd x\dd y\nonumber\\ 
	& \leq   C \big\langle A_{s_1,p}(\underline{u}_{\lambda}) + A_{s_2,q}(\underline{u}_{\lambda}),\phi \big\rangle \nonumber\\
	&= C\int_{\Omega_{k}} h_{\lambda}(\underline{u}_{\lambda})(\underline{u}_{\lambda}-k)\dd x\\
	&= C\int_{\Omega_{k}} \big[\lambda c_0 \underline{u}_{\lambda}^{\tau - 1} - \lambda c_2\underline{u}_{\lambda}^{\theta-1}\big](\underline{u}_{\lambda}-k)\dd x \nonumber\\
	&\leq C\int_{\Omega_{k}} \lambda c_0 \underline{u}_{\lambda}^{\tau - 1} (\underline{u}_{\lambda}-k)\dd x. \nonumber
\end{align}
Since $1<\tau<p$, for $k>1$ in $\Omega_k$ we have 
\[
\underline{u}_{\lambda}^{\tau-1}(\underline{u}_{\lambda}-k)\leq \underline{u}_{\lambda}^{p-1}(\underline{u}_{\lambda}-k)\leq2^{p-1}(\underline{u}_{\lambda}-k)^p+2^{p-1}k^{p-1}(\underline{u}_{\lambda}-k)
\]
and thus,
\begin{equation}\label{eq 4.20}
	\int_{\Omega}\underline{u}_{\lambda}^{\tau-1}(\underline{u}_{\lambda}-k)\dd x\leq 2^{p-1}\int_{\Omega}(\underline{u}_{\lambda}-k)^p \dd x+2^{p-1}k^{p-1}\int_{\Omega_k}(\underline{u}_{\lambda}-k)\dd x.
\end{equation}

Applying  Hölder's inequality, we obtain
\begin{equation}\label{eq 4.22}
	\int_{\Omega_k}(\underline{u}_{\lambda}-k)^p\dd x\leq |\Omega_k|^{\frac{p^*_s-p}{p^*_s}}\left( \int_{\Omega_k}(\underline{u}_{\lambda}-k)^{p^*_s} \dd x\right)^{\frac{p}{p^*_s}}.
\end{equation}

So, using the inequalities \eqref{eq 4.20} and \eqref{eq 4.22} in \eqref{eq1}, we have
$$\int_{\Omega_k}(\underline{u}_{\lambda} -k)^p\dd x \leq C_0\left|\Omega_k\right|^{\frac{p^*_s-p}{p^*_s}}\left[2^{p-1}\int_{\Omega_k}(\underline{u}_{\lambda}-k)^p\dd x+2^{p-1}k^{p-1}\int_{\Omega_k}(\underline{u}_{\lambda}-k)\dd x \right].$$
Thus, we obtain
\begin{equation*}
	\left[1- 2^{p-1}C_0|\Omega_k|^{\frac{p^*_s-p}{p^*_s}}\right] \int_{\Omega_k}(\underline{u}_{\lambda} - k)^p\dd x\leq 2^{p-1}k^{p-1}\left|\Omega_k\right|^{\frac{(p^*_s-p)}{p^*_s}}\int_{\Omega_k}(\underline{u}_{\lambda} - k)\dd x.
\end{equation*}

If $k\to \infty$, then $|\Omega_k|\to 0$. Therefore, there exists $k_0> 1$ such that
$$1-2^{p-1}C_0|\Omega_k|^{\frac{p^*_s-p}{p^*_s}}\geq\frac{1}{2}\quad\textrm{if}~~k\geq k_0>1.$$

Thus, for such $k$, we conclude that
\begin{equation}\label{eq 4.24}
	\frac{1}{2}\int_{\Omega_k}(\underline{u}_{\lambda}-k)^p\dd x \leq  2^{p-1}k^{p-1}C_0|\Omega_k|^{\frac{p^*_s-p}{p^*_s}}\int_{A_k}(\underline{u}_{\lambda}-k)\dd x.
\end{equation}

Hölder's inequality and \eqref{eq 4.24} yield
\begin{equation*}
	\left(\int_{\Omega_k}(\underline{u}_{\lambda}-k)\dd x \right)^{p} \leq  |\Omega_k|^{p-1}\int_{\Omega_k}(\underline{u}_{\lambda}-k)^p\dd x\leq  |\Omega_k|^{p-1}2^{p-1}k^{p-1}C_0|\Omega_k|^{\frac{p^*_s-p}{p^*_s}}\int_{A_k}(\underline{u}_{\lambda}-k)\dd x.
\end{equation*}

Thus,
\begin{equation}\label{eq 4.26}
	\int_{\Omega_k}(u-k)\dd x\leq 2\tilde{C} k|\Omega_k|^{1+\epsilon},\quad\forall\, k\geq k_0,
\end{equation}
where $\epsilon=\displaystyle\frac{p^*_s-p}{p^*_s{(p-1)}}>0$ and $\tilde{C} >0$. 

The same arguments used in \cite{Gilberto} assures us that $\underline{u}_{\lambda} \in L^{\infty}(\Omega)$. Then the nonlinear regularity theory, see \cite{Gio} says that $\underline{u}_{\lambda} \in \text{int}(C_{s_1}^0(\Omega))_+$.

\qed

\begin{propo}\label{prop2}
For every $\lambda>0$, the problem \eqref{P2} admits a unique positive solution $\underline{u}_\lambda \in$ int$(C^0_{s_1}(\Omega)_+)$ and $\underline{u}_{\lambda} \to 0$ in $C_{s_1}^{0}(\overline{\Omega})$ as $\lambda \to 0^+$.
\end{propo}
\n {\bf Proof.} \textbf{Existence:} Note that, the solutions of the problem \eqref{P2} are critical points of the functional $\tilde{I}_{\lambda}: W_{0}^{s_1,p}(\Omega) \to W_{0}^{s_2,q}(\Omega)$ given by
\begin{equation}
\tilde{I}_\lambda(u) = \frac{1}{p} \left[ u \right]_{s_1,p}^p + \frac{1}{q} \left[ u \right]_{s_2,p}^{q} - \frac{\lambda c_0}{\tau} \Vert u^+ \Vert_{\tau}^{\tau} + \frac{\lambda c_2}{\theta} \Vert u^+ \Vert_\theta^\theta, \ \text{for all} \ u \in W_0^{s_1,p}(\Omega)
\end{equation}
where $\Vert.\Vert_{t}$ denote the norm in space $L^{t}(\Omega)$.

Since $1<\tau<q<p<\theta$, then $\tilde{I}_\lambda(tu) \to \infty$ as $t \to \infty$, is that, $J_\lambda$ is coercive. Also using the Sobolev embedding theorem, we see that $\tilde{I}_\lambda$ is sequentially weakly lower semicontinuous. So, by the Weierstrass-Tonelli theorem, we can find $\underline{u}_\lambda \in W_0^{s_1,p}(\Omega)$ such that
\[
\tilde{I}_\lambda(\underline{u}_\lambda) = \min\bigg\{J_\lambda(u); \ \ u \in W_0^{s_1,p}(\Omega)\bigg\}.
\]

Now notice that $1<\tau<q<p<\theta$ and $u\in$ int$(C^0_{s_1}(\Omega)_+)$ results
\begin{equation}\label{min}
\tilde{I}_{\lambda}(tu) < 0 \ \ \text{for} \ \ t \in (0,1) \ \ \text{small enough}
\end{equation}
thus $\tilde{I}_{\lambda}(\underline{u}_{\lambda}) < 0 = \tilde{I}_{\lambda}(0)$ and therefore $\underline{u}_{\lambda} \neq 0$.

Using the \eqref{min} we have,
\[
\tilde{I}_{\lambda}'(\underline{u}_{\lambda}) = 0
\]
and consequently
\begin{equation}\label{test}
\bigg\langle A_{s_1,p}(\underline{u}_{\lambda}) + A_{s_2,q}(\underline{u}_{\lambda}),\varphi \bigg\rangle = \lambda \int_{\Omega} c_0(\underline{u}^{+}_{\lambda})^{\tau -1}\varphi \dd x - \lambda \int_{\Omega} c_2(\underline{u}^{+}_{\lambda})^{\theta-1}\varphi \dd x, \ \ \text{for all} \ \ \varphi \in W_{0}^{s_1,p}(\Omega).
\end{equation}

Choosing $\varphi = \underline{u}_{\lambda}^{-} \in W_{0}^{s_1,p}(\Omega)$ results
\[
\left[\underline{u}_{\lambda}^{-} \right]_{s_1,p}^{p} + \left[ \underline{u}_{\lambda}^{-} \right]_{s_2,q} \leq \bigg\langle A_{s_1,p}(\underline{u}_{\lambda}) + A_{s_2,q}(\underline{u}_{\lambda}),\underline{u}_{\lambda}^{-} \bigg\rangle = \lambda \int_{\Omega} c_0(\underline{u}^{+}_{\lambda})^{\tau -1}\underline{u}_{\lambda}^{-} \dd x - \lambda \int_{\Omega} c_2(\underline{u}^{+}_{\lambda})^{\theta-1}\underline{u}_{\lambda}^{-} \dd x = 0
\]
and therefore $\left[ \underline{u}_{\lambda}^{-} \right]_{s_1,p}^{p} = 0$, is that, $ \underline{u}_{\lambda} \geq 0$ and $\underline{u}_{\lambda} \neq 0$.




\textbf{Uniqueness:} To show the uniqueness of the solution, we will use arguments similar to those used in \cite{IM}. Let's use the following discrete Picone's inequality from \cite{Picone}
\begin{equation}\label{Picone}
J_r(a-b)\left(\frac{c^r}{a^{r-1}} - \frac{d^{r}}{b^{r-1}}\right) \leq \vert c-d \vert^{r}, \ \ \text{for all} \ a,b \in \mathbb{R}^{*}_+, c,d \in \mathbb{R}^{+}.
\end{equation}

Let $\underline{u}_{\lambda},\underline{v}_{\lambda} \in W_{0}^{s_1,p}(\Omega)$ positive solutions of the problem \eqref{P2}. As above, we show that $\underline{u}_\lambda,\underline{v}_\lambda \in$ int$(C^0_{s_1}(\Omega)_+)$. Thus, using the same arguments as Lemma 2.4 of \cite{IM} we have,
\[
\frac{\underline{u}_\lambda^p}{\underline{v}_\lambda^{p-1}} \in W_0^{s_1,p}(\Omega).
\]

Consider $w_{\lambda} = (\underline{u}_{\lambda}^{p} - \underline{v}_{\lambda}^{p})^{+}$, thus,
\[
\frac{w_{\lambda}}{\underline{v}_{\lambda}^{p - 1}} = \left(\frac{\underline{u}_{\lambda}^{p}}{\underline{v}_{\lambda}^{p-1}} - \underline{v}_{\lambda}\right)^{+} \in W_{0}^{s_1,p}(\Omega) \ \text{and} \ \frac{w_{\lambda}}{\underline{u}_{\lambda}^{p - 1}} = \left(\underline{u}_{\lambda} -\frac{\underline{v}_{\lambda}^{p}}{\underline{u}_{\lambda}^{p-1}}\right)^{+} \in W_{0}^{s_1,p}(\Omega).
\]

We denote by $g_{\lambda}(t) = \lambda c_0 t^{\tau - p} - \lambda c_2t^{\theta-p}$. Thus, $g$ is strictly dcreasing in $\mathbb{R}^{+}_0$.

Testing \eqref{test} with $\frac{w_{\lambda}}{\underline{u}_{\lambda}^{p - 1}}$ we have
\begin{align*}
\bigg\langle A_{s_1,p}(\underline{u}_{\lambda}) + A_{s_2,q}(\underline{u}_{\lambda}), \frac{w_{\lambda}}{\underline{u}_{\lambda}^{p - 1}} \bigg\rangle &= \lambda \int_{\Omega} c_0\underline{u}_{\lambda}^{\tau -1}\frac{w_{\lambda}}{\underline{u}_{\lambda}^{p - 1}} \dd x - \lambda \int_{\Omega} c_2\underline{u}_{\lambda}^{\theta-1}\frac{w_{\lambda}}{\underline{u}_{\lambda}^{p - 1}} \dd x \\
&= \lambda \int_{\Omega} c_0\underline{u}_{\lambda}^{\tau-p} w_{\lambda} \dd x - \lambda \int_{\Omega} c_2\underline{u}_{\lambda}^{\theta-p}w_{\lambda}\dd x \\
&= \int_{\{\underline{u}_{\lambda} > \underline{v}_{\lambda}\}} g_{\lambda}(\underline{u}_{\lambda}) (\underline{u}_{\lambda}^p - \underline{v}_{\lambda}^p) \dd x
\end{align*}
and testing \eqref{test} with 
$\frac{w_{\lambda}}{\underline{v}_{\lambda}^{p - 1}}$ we have
\begin{align*}
\bigg\langle A_{s_1,p}(\underline{v}_{\lambda}) + A_{s_2,q}(\underline{v}_{\lambda}), \frac{w_{\lambda}}{\underline{v}_{\lambda}^{p - 1}} \bigg\rangle &= \lambda \int_{\Omega} c_0\underline{v}_{\lambda}^{\tau -1}\frac{w_{\lambda}}{\underline{v}_{\lambda}^{p - 1}} \dd x - \lambda \int_{\Omega} c_2\underline{v}_{\lambda}^{\theta-1}\frac{w_{\lambda}}{\underline{v}_{\lambda}^{p - 1}} \dd x \\
&= \lambda \int_{\Omega} c_0 \underline{v}^{\tau - p} w_{\lambda} \dd x - \lambda \int_{\Omega} c_2\underline{v}_{\lambda}^{\theta-\tau}w_{\lambda}\dd x\\
&= \int_{\{\underline{u}_{\lambda} > \underline{v}_{\lambda}\}} g_{\lambda}(\underline{v}_{\lambda}) (\underline{u}_{\lambda}^p - \underline{v}_{\lambda}^p) \dd x
\end{align*}

Thus,
\begin{align*}
\bigg\langle A_{s_1,p}(\underline{u}_{\lambda}) + A_{s_2,q}(\underline{u}_{\lambda}), \frac{w_{\lambda}}{\underline{u}_{\lambda}^{p - 1}} \bigg\rangle - \bigg\langle A_{s_1,p}(\underline{v}_{\lambda}) + A_{s_2,q}(\underline{v}_{\lambda}), \frac{w_{\lambda}}{\underline{v}_{\lambda}^{p - 1}} \bigg\rangle = \int_{\{\underline{u}_{\lambda} > \underline{v}_{\lambda}\}} \left[g_{\lambda}(\underline{u_\lambda}) - g_{\lambda}(\underline{v}_{\lambda})\right] (\underline{u}_{\lambda}^p - \underline{v}_{\lambda}^p) \dd x.
\end{align*}

Note that, using the discrete Picone's inequality \eqref{Picone}, see (Proposition 3.1, \cite{IM}) we have
\[
j_p(u(x) - u(y))\left(\frac{w_{\lambda}(x)}{\underline{u}_{\lambda}(x)^{p-1}} - \frac{w_{\lambda}(y)}{\underline{u}_{\lambda}(y)^{p-1}} \right) \geq j_p(v(x) - v(y))\left(\frac{w_{\lambda}(x)}{\underline{v}_{\lambda}(x)^{p-1}} - \frac{w_{\lambda}(y)}{\underline{v}_{\lambda}(y)^{p-1}} \right)
\]
and thus,
\[
\bigg\langle A_{s_1,p}(\underline{u}_{\lambda}) + A_{s_2,q}(\underline{u}_{\lambda}), \frac{w_{\lambda}}{\underline{u}_{\lambda}^{p - 1}} \bigg\rangle \geq \bigg\langle A_{s_1,p}(\underline{v}_{\lambda}) + A_{s_2,q}(\underline{v}_{\lambda}), \frac{w_{\lambda}}{\underline{v}_{\lambda}^{p - 1}} \bigg\rangle.
\]

Therefore, since $g_{\lambda}$ is strictly decreasing in $\mathbb{R}^{+}_0$ results
\[
0 \leq \int_{\{\underline{u}_{\lambda} > \underline{v}_{\lambda}\}} \left[g_{\lambda}(\underline{u_\lambda}) - g_{\lambda}(\underline{v}_{\lambda})\right] (\underline{u}_{\lambda}^p - \underline{v}_{\lambda}^p) \dd x \leq 0
\]
so we deduce that $\{\underline{u}_{\lambda} > \underline{v}_{\lambda}\}$ has null measure, is that, $\underline{u}_{\lambda} \leq \underline{v}_{\lambda}$ in $\Omega$. Similarly, using the function test $w_\lambda = (\underline{v}_{\lambda}^p - \underline{u}_{\lambda}^p)^+$ we see that $\underline{u}_{\lambda} \geq \underline{v}_{\lambda}$ in $\Omega$, and thus $\underline{u}_{\lambda} = \underline{v}_{\lambda}$.

Moreover, we have
\begin{align*}
\left[ \underline{u}_{\lambda} \right]_{s_1,p}^{p} &\leq \left[\underline{u}_{\lambda} \right]_{s_1,p}^{p} + \left[\underline{u}_{\lambda} \right]_{s_2,q}^{q} \\
&= \lambda c_0 \Vert \underline{u}_{\lambda} \Vert_{\tau}^{\tau} - \lambda c_2 \Vert \underline{u}_{\lambda} \Vert_{\theta}^{\theta}\\
&\leq \lambda c_0\Vert \underline{u}_{\lambda} \Vert_{\tau}^{\tau}\\
&\leq \lambda \hat{c}_{0} \left[\underline{u}_{\lambda} \right]_{s_1,p}^{\tau},
\end{align*}
for some $\hat{c}_0 >0$. Thus,
\[
\left[\underline{u}_{\lambda} \right]_{s_1,p}^{p-\tau} \leq \lambda\hat{c}_{0} 
\]
and therefore, $\underline{u}_{\lambda} \to 0$ in $W^{s_1,p}_0(\Omega)$ as $\lambda \to 0^+$. Using the nonlinear regularity theorem, see\cite{Gio}, results that
\[
\underline{u}_{\lambda} \to 0 \ \ \text{in} \ \ C_{s_1}^{0}(\overline{\Omega}) \ \ \text{as} \ \ \lambda \to 0^+.
\]
\qed

We consider another auxiliary problem,
\begin{equation} \label{P3}
\left\{
\begin{array}{llll}
(-\Delta_p)^{s_1}u + (-\Delta_q)^{s_2}u = \lambda \underline{u}_{\lambda}^{-\eta} + 1 & {\rm in} \ \ \Omega,\\
u= 0 & {\rm in} \ \R^N\setminus\Omega,\\
u> 0 & {\rm in} \ \Omega
\end{array}
\right.
\end{equation}
with $\lambda > 0$, $0<\eta<1$ and $1 < q<p$.
\begin{propo}\label{prop3}
For every $\lambda > 0$, there exists a unique solution $\overline{u}_{\lambda} \in \text{int}\left[(C_{s_1}^{0}(\overline{\Omega}))_+\right]$ of the problem \eqref{P3} and a $\lambda_0>0$ such that, for all $0 < \lambda \leq \lambda_0$ it holds
\[
\underline{u}_{\lambda} \leq \overline{u}_{\lambda}.
\]
\end{propo}
\n {\bf Proof.} Note that, the Lemma 14.16 of Gilbarg-Trundiger \cite{GT} says that $d^{s_1}_{\Omega} \in C^{2}(\Omega_{\delta_0})$, where $\Omega_{\delta_0} = \{x \in \Omega; d^{s_1}_{\Omega}(x) < \delta_0\}$. Thus, $d^{s_1}_{\Omega} \in \text{int}\left[(C_{s_1}^0(\Omega))_+\right]$ and so by Proposition 4.1.22 of \cite{Radulesco}, there exists $c_3 = c_3(\underline{u}_{\lambda}) >0$ and $c_4 = c_4(\underline{u}_{\lambda}) > 0$ such that,
\begin{equation}\label{est_dist_f}
c_3 d^{s_1}_{\Omega} \leq \underline{u}_{\lambda} \leq c_4d^{s_1}_{\Omega}.
\end{equation}

Since due to (\ref{est_dist_f}), $\lambda \underline{u}_{\lambda}^{-\eta} + 1 \in L^1(\Omega)$. The existence of a weak solution of (\ref{P3}) follows from direct minimization in $W_0^{s_1, p}(\Omega)$ of the functional
\[
\frac{1}{p} \left[ u \right]_{s_1,p}^p + \frac{1}{q} \left[ u \right]_{s_2,p}^{q} - \int_{\Omega}(\lambda \underline{u}_{\lambda}^{-\eta} + 1)udx,
\]
whereas the uniqueness comes from, for instance, the comparison principle for the Dirichlet fractional $(p,q)$-Laplacian, Propossition \ref{PC}. Using the maximum principle, \cite{Gio}, the solution $\overline{u}_{\lambda} \in \text{int}\left[(C_{s_1}^{0}(\overline{\Omega})_+)\right]$.

For show the existence of $\lambda_0 >0$ such that $\underline{u}_{\lambda} \leq \overline{u}_{\lambda}$ for all $0< \lambda \leq \lambda_0$, acting on \eqref{P3} with $\overline{u}_{\lambda}$ and obtain
\begin{align*}
\left[ \overline{u}_{\lambda} \right]_{s_1,p}^p &\leq \left[ \overline{u}_{\lambda} \right]_{s_1,p}^p + \left[ \overline{u}_{\lambda} \right]_{s_2,q}^q\\
&= \lambda \displaystyle\int_{\Omega} \underline{u}_{\lambda}^{-\eta}.\overline{u}_{\lambda} \dd x + \displaystyle\int_{\Omega} \overline{u}_{\lambda} \dd x \\
&= \lambda \displaystyle\int_{\Omega} \underline{u}_{\lambda}^{1-\eta}.\frac{\overline{u}_{\lambda}}{\underline{u}_{\lambda}} \dd x + \displaystyle\int_{\Omega} \overline{u}_{\lambda} \dd x \\ 
&\leq \lambda c_5 \int_{\Omega} \frac{\overline{u}_{\lambda}}{d^{s_1}_{\Omega}} \dd x + \vert \Omega \vert^{\frac{p-1}{p}} \left(\int_{\Omega} \overline{u}_{\lambda}^{p} \dd x\right)^{\frac{1}{p}} \ \ (\text{Holder inequality})\\
&\leq \left(\lambda c_5 + \frac{1}{\lambda_1(p)}\right)\vert \Omega\vert^{\frac{p-1}{p}} \left[ \overline{u}_{\lambda} \right]_{s_1,p} (\text{Hardy's inequality and first eigenvalue}).
\end{align*} 

So, we have $\{\overline{u}_{\lambda}\}_{\lambda \in (0,1]}$ is uniformly bounded in $W_{0}^{s_1,p}(\Omega)$. Using arguments similar to the Lemma \ref{Lema1}, (see also Ladyzhenskaya-Ural'tseva \cite{Lady} Theorem 7.1) results 
\[
\{\overline{u}_{\lambda}\}_{\lambda \in (0,1]} \subset L^{\infty}(\Omega) \ \ \text{is uniformly bounded in} \ \ \lambda.
\]

The condition $\textbf{H} \ (i)$ implies that there exists $\lambda_0>0$ such that,
\[
\lambda f_k(x,\overline{u}_{\lambda}) \leq \lambda \Vert a \Vert (1 + \Vert \overline{u}_{\lambda} \Vert^{\theta-1}) \leq 1 \ \ \text{for all} \ \ \lambda \in (0,\lambda_0] \ \ \text{and} \ \ x \ \ \text{a. a. in} \ \Omega.
\]

For each $\lambda \in (0,\lambda_0]$ consider the Carathéodory function 
\[
\kappa_{\lambda}(x,t) = \left\{
\begin{array}{llll}
\lambda[c_0(t^{+})^{\tau -1} - c_2(t^{+})^{\theta-1}] & {\rm if} \ \ t \leq \overline{u}_{\lambda}(x),\\
\lambda[c_0\overline{u}_{\lambda}(x)^{\tau -1} - c_2\overline{u}_{\lambda}(x)^{\theta-1}] & {\rm if} \ \  \overline{u}_{\lambda}(x) < t.
\end{array}
\right.
\]

Let $\Psi_{\lambda}: W_{0}^{s_1,p} \to \mathbb{R}$ the $C^1$-functional defined by
\[
\Psi_{\lambda}(u) = \frac{1}{p} \left[ u \right]_{s_1,p}^{p} + \frac{1}{q} \left[ u \right]_{s_2,p}^{q} - \int_{\Omega} K_{\lambda}(x,u) \dd x, \ \ \text{for all} \ \ u \in W_{0}^{s_1,p}(\Omega)
\]
where $K_{\lambda}(x,t) = \displaystyle\int_{0}^{t} \kappa_\lambda(x,s) \dd s$.

Note that, $\Psi_{\lambda}$ is coercive and sequentially wekly lower semicontinuous. So, there exists $\tilde{u}_{\lambda} \in W_0^{s_1,p}(\Omega)$ such that
\[
\Psi_{\lambda}(\tilde{u}_{\lambda}) = \text{min}\left[\Psi_{\lambda}(u); \ u \in W_0^{s_1,p}(\Omega)\right].
\]

Since $1<\tau<q<p<\theta$ results
\begin{equation}\label{min1}
\Psi_{\lambda}(tu) < 0 \ \ \text{for} \ \ t \in (0,1) \ \ \text{small enough}
\end{equation}
thus $\Psi_{\lambda}(\underline{u}_{\lambda}) < 0 = \Psi_{\lambda}(0)$ and therefore $\underline{u}_{\lambda} \neq 0$.

Using the \eqref{min1} we have,
\[
\Psi_{\lambda}'(\tilde{u}_{\lambda}) = 0
\]
and consequently
\begin{equation}\label{test}
\bigg\langle A_{s_1,p}(\tilde{u}_{\lambda}) + A_{s_2,q}(\tilde{u}_{\lambda}),\varphi \bigg\rangle = \int_{\Omega}\kappa_{\lambda}(x,\tilde{u}_{\lambda}) \varphi \dd x, \ \ \text{for all} \ \ \varphi \in W_{0}^{s_1,p}(\Omega).
\end{equation}

Choosing $\varphi = -\tilde{u}_{\lambda} \in W_{0}^{s_1,p}(\Omega)$, we see that $\tilde{u}_{\lambda} \geq 0$ and $\tilde{u}_{\lambda} \neq 0$. Taking $\varphi = (\tilde{u}_{\lambda} - \overline{u}_{\lambda})^+ \in W_0^{s_1,p}(\Omega)$ we find, 

From \eqref{fcond}, we have that there exits $c_0>0$ and $c_2>0$ such that $f_k(x,t) \geq c_0t^{\tau-1} - c_2t^{\theta-1}$ and so
\begin{align*}
\bigg\langle A_{s_1,p}(\tilde{u}_{\lambda}) + A_{s_2,q}(\tilde{u}_{\lambda}),(\tilde{u}_{\lambda} - \overline{u}_{\lambda})^+ \bigg\rangle &= \int_{\Omega} \kappa_{\lambda}(x,\tilde{u}_{\lambda})(\tilde{u}_{\lambda} - \overline{u}_{\lambda})^+ \dd x\\
&= \int_{\Omega} \lambda[c_0\overline{u}_{\lambda}^{\tau -1} - c_2\overline{u}_{\lambda}^{\theta-1}](\tilde{u}_{\lambda} - \overline{u}_{\lambda})^+ \dd x\\
&\leq \int_{\Omega} \lambda f_k(x,\overline{u}_{\lambda})(\tilde{u}_{\lambda} - \overline{u}_{\lambda})^+ \dd x\\
&\leq \int_{\Omega} [\lambda \underline{u}_{\lambda}^{-\eta} + 1](\tilde{u}_{\lambda} - \overline{u}_{\lambda})^+ \dd x    \ \ \ (\text{for all} \ \ 0< \lambda \leq \lambda_0)\\
&= \bigg\langle A_{s_1,p}(\overline{u}_{\lambda}) + A_{s_2,q}(\overline{u}_{\lambda}),(\tilde{u}_{\lambda} - \overline{u}_{\lambda})^+ \bigg\rangle
\end{align*}
and so, by Proposition \ref{PC} $\tilde{u}_{\lambda} \leq \overline{u}_{\lambda}$. Moreover, note that,
\[
\Psi_{\lambda}(u) = \tilde{I}_{\lambda}(u), \ \ \text{for all} \ \ u \in [0,\overline{u}_{\lambda}],
\]
thus
\begin{align*}
\tilde{I}_\lambda(\tilde{u}_{\lambda}) &= \Psi_{\lambda}(\tilde{u}_\lambda) = \text{min}\left[\Psi_{\lambda}(u); \ u \in W_0^{s_1,p}(\Omega)\right]\\
&= \text{min}\bigg\{\Psi_{\lambda}(u); \ u \in [0,\overline{u}_{\lambda}]\bigg\}\\
&= \text{min}\bigg\{\tilde{I}_{\lambda}(u); \ u \in [0,\overline{u}_{\lambda}]\bigg\}\\
&= \tilde{I}_{\lambda}(\underline{u}_{\lambda}).
\end{align*}

By Proposition \ref{prop2} we have $\tilde{u}_{\lambda} = \underline{u}_{\lambda}$ and therefore $\underline{u}_{\lambda} \leq \overline{u}_{\lambda}$ for all $0 < \lambda \leq \lambda_0$.
\qed

\section{Existence of positive solution for \ref{P1}}

We consider the set 
\[
\mathcal{L} = \bigg\{\lambda>0; \ \ \text{problem \ref{P1} admits a positive solution}\bigg\}
\]
and the set $S_{\lambda}$ of the positive solutions to the problem \ref{P1}.
\begin{propo}\label{propex}
Assume the hypoteses $(\textbf{H}_k)$ hold, then
\begin{enumerate}
\item[i)] $\mathcal{L} \neq \varnothing$;
\item[ii)] If $\lambda \in \mathcal{L}$, then $\underline{u}_{\lambda} \leq u$ for all $u \in S_{\lambda}$ and $S_{\lambda} \subseteq \text{int}[(C_{s_1}^0(\Omega))_+]$.
\end{enumerate}
\end{propo}
\n {\bf Proof.} Let $\lambda_0 > 0$ given in the Proposition \ref{prop2}, so for $\lambda \in (0,\lambda_0]$ we have
\begin{equation}\label{eq3}
\underline{u}_{\lambda} \leq \overline{u}_{\lambda} \ \ \text{and} \ \ \lambda f(x,\overline{u}_{\lambda}) \leq 1 \ \ \text{for a. a.} \ \ x \in \Omega.
\end{equation}

We consider the function 
\[
g_{\lambda}(x,t) = \left\{
\begin{array}{llll}
\lambda[ \underline{u}_{\lambda}^{-\eta} + f_k(x, \underline{u}_{\lambda})] & {\rm if} \ \ t < \underline{u}_{\lambda}(x),\\
\lambda[t^{-\eta} + f_k(x, t)] & {\rm if} \ \  \underline{u}_{\lambda}(x) \leq t \leq \overline{u}_{\lambda}(x),\\
\lambda[\overline{u}_{\lambda}^{-\eta} + f_k(x,\overline{u}_{\lambda})] & {\rm if} \ \ \overline{u}_{\lambda}(x)< t,
\end{array}
\right.
\]
and the functional $\Phi_\lambda: W_0^{s_1,p}(\Omega) \to \mathbb{R}$ defined by
\[
\Phi_{\lambda}(u) = \frac{1}{p} \left[ u \right]_{s_1,p}^{p} + \frac{1}{q} \left[ u \right]_{s_2,p}^{q} - \int_{\Omega} G_{\lambda}(x,u) \dd x, \ \ \text{for all} \ \ u \in W_{0}^{s_1,p}(\Omega)
\]
where $G(x,t) = \displaystyle\int_{0}^{t} g_{\lambda}(x,s) \dd s$.

By Proposition $3$ of \cite{Papa} we have $\Phi_\lambda \in C^{1}(W_0^{s_1,p}(\Omega), \mathbb{R})$. Moreover, using the hipoteses $(\textbf{H})$ we have, $\Phi_{\lambda}$ is coercive and sequently weakly lower semicontinuous. Thus, there exists $u_{\lambda}:= u_{k,\lambda} \in W_0^{s_1,p}(\Omega)$ such that,
\[
\Phi_{\lambda}(u_{\lambda}) = \text{min}\bigg[\Phi_{\lambda}(u); \ u \in W_0^{s_1,p}(\Omega)\bigg].
\]
Thus, $\Phi_\lambda'(u_\lambda) = 0$, that is,
\begin{equation}\label{test1}
\bigg\langle A_{s_1,p}(u_{\lambda}) + A_{s_2,q}(u_{\lambda}),\varphi \bigg\rangle = \int_{\Omega}g_{\lambda}(x,u_{\lambda}) \varphi \dd x, \ \ \text{for all} \ \ \varphi \in W_{0}^{s_1,p}(\Omega).
\end{equation}

Testing the equation \eqref{test1} with $\varphi = (u_\lambda-\overline{u}_{\lambda})^+ \in W_0^{s_1,p}(\Omega)$ and using the inequality \eqref{eq3}, we find 
\begin{align*}
\bigg\langle A_{s_1,p}(u_{\lambda}) + A_{s_2,q}(u_{\lambda}),(u_{\lambda} - \overline{u}_{\lambda})^+ \bigg\rangle &= \int_{\Omega} g_{\lambda}(x,u_{\lambda})(u_{\lambda} - \overline{u}_{\lambda})^+ \dd x\\
&= \int_{\Omega} \lambda[\overline{u}_{\lambda}^{-\eta} + f_k(x,\overline{u}_{\lambda})](u_{\lambda} - \overline{u}_{\lambda})^+ \dd x\\
&\leq \int_{\Omega} [\lambda \underline{u}_{\lambda}^{-\eta} + 1](u_{\lambda} - \overline{u}_{\lambda})^+ \dd x    \ \ \ (\text{for all} \ \ 0< \lambda \leq \lambda_0)\\
&= \bigg\langle A_{s_1,p}(\overline{u}_{\lambda}) + A_{s_2,q}(\overline{u}_{\lambda}),(u_{\lambda} - \overline{u}_{\lambda})^+ \bigg\rangle
\end{align*}
and so, by Proposition \ref{PC} $u_{\lambda} \leq \overline{u}_{\lambda}$.

Analogously, testing \eqref{test1} with the function $\varphi = (\underline{u}_{\lambda} - u_\lambda)^+ \in W_0^{s_1,p}(\Omega)$ and using \eqref{fcond}, we have,
\begin{align*}
\bigg\langle A_{s_1,p}(u_{\lambda}) + A_{s_2,q}(u_{\lambda}),(\underline{u}_{\lambda} - u_{\lambda})^+ \bigg\rangle &= \int_{\Omega} g_{\lambda}(x,u_{\lambda})(\underline{u}_{\lambda} - u_{\lambda})^+ \dd x\\
&= \int_{\Omega} \lambda[\underline{u}_{\lambda}^{-\eta} + f_k(x,\underline{u}_{\lambda})](\underline{u}_{\lambda} - u_{\lambda})^+ \dd x\\
&\geq \int_{\Omega} \lambda[c_0\underline{u}_\lambda^{\tau -1} - c_2\underline{u}^{\theta-1}](\underline{u}_{\lambda} - u_{\lambda})^+ \dd x    \ \ \ (\text{for all} \ \ 0< \lambda \leq \lambda_0)\\
&= \bigg\langle A_{s_1,p}(\underline{u}_{\lambda}) + A_{s_2,q}(\underline{u}_{\lambda}),(\underline{u}_{\lambda} - u_{\lambda})^+ \bigg\rangle
\end{align*}
and so, by Proposition \ref{PC} we have $u_{\lambda} \leq \overline{u}_{\lambda}$.

Therefore, 
\[
u_\lambda \in [\underline{u}_\lambda,\overline{u}_\lambda] \Rightarrow u_\lambda \in S_\lambda \Rightarrow (0,\lambda_0] \subseteq \mathcal{L}.
\]

For item (ii), it is sufficient to argue as in the Proposition \ref{prop2}, replacing $\overline{u}_\lambda$ with $u \in S_\lambda$, we show that $\underline{u}_\lambda \leq u$ for all $u \in S_\lambda$. For show that $S_\lambda \subseteq \text{int}[(C_{s_1}^0(\Omega))_+]$ we use the maximum principle, see \cite{Gio}.
\qed
\begin{propo}
If hypotheses $(\textbf{H}_k)$ hold, $\lambda \in \mathcal{L}$ and $\mu \in(0,\lambda)$, then $\mu \in \mathcal{L}$.
\end{propo}
\n {\bf Proof.} Let $\lambda \in \mathcal{L}$, so we can find $u_{\lambda} \in S_{\lambda} \subseteq \text{int}[(C_{s_1}^{0}(\overline{\Omega}))_+]$. Consider the Dirichlet problem,
\begin{equation} \label{P4}
	\left\{
	\begin{array}{llll}
		(-\Delta_p)^{s_1}u + (-\Delta_q)^{s_2}u = \vartheta c_0u(x)^{\tau -1} - \lambda c_2u^{\theta -1} & {\rm in} \ \ \Omega,\\
		u= 0 & {\rm in} \ \R^N\setminus\Omega,\\
		u> 0 & {\rm in} \ \Omega
	\end{array}
	\right.
\end{equation}
with $0 < \vartheta<\lambda$ and $1 < \tau < q<p<\theta$. As we did in the proposition, we can find a unique solution $\tilde{u}_{\vartheta} \in \text{int}[(C_{s_1}^{0}(\overline{\Omega}))_+]$ to the problem \eqref{P4} and, in addition, we can show that $\tilde{u}_{\vartheta}^{-\eta} \in L^{1}(\Omega)$. Since, for all $0<\vartheta_1 < \vartheta_2 \leq \lambda$, we have $\vartheta_1 c_0u(x)^{\tau -1} - \lambda c_2u^{\theta -1} \leq \vartheta_2 c_0u(x)^{\tau -1} - \lambda c_2u^{\theta -1}$, by comparison principle results that $\tilde{u}_{\vartheta_1} \leq \tilde{u}_{\vartheta_2}$. Note that, by Proposition \ref{prop3} $\tilde{u}_{\lambda} = \underline{u}_{\lambda}$, so
\[
\tilde{u}_{\mu} \leq \underline{u}_{\lambda} \leq u_{\lambda}.
\]

Define the Caracthéodoryfunction,
\[
\gamma_{\lambda}(x,t) = \left\{
\begin{array}{llll}
	\mu[ \tilde{u}_{\mu}^{-\eta} + f_k(x, \tilde{u}_{\mu})] & {\rm if} \ \ t < \tilde{u}_{\mu}(x),\\
	\mu[t^{-\eta} + f_k(x, t)] & {\rm if} \ \  \tilde{u}_{\mu}(x) \leq t \leq \tilde{u}_{\mu}(x),\\
	\mu[\tilde{u}_{\mu}^{-\eta} + f_k(x,\tilde{u}_{\mu})] & {\rm if} \ \ \tilde{u}_{\mu}(x)< t,
\end{array}
\right.
\]

Let $\Upsilon_{\lambda}: W_{0}^{s_1,p}(\Omega) \to \mathbb{R}$ the $C^1$-functional defined by
\[
\Upsilon_{\lambda}(u) = \frac{1}{p} \left[ u \right]_{s_1,p}^{p} + \frac{1}{q} \left[ u \right]_{s_2,p}^{q} - \int_{\Omega} \Gamma_{\mu}(x,u) \dd x, \ \ \text{for all} \ \ u \in W_{0}^{s_1,p}(\Omega)
\]
where $\Gamma_{\lambda}(x,t) = \displaystyle\int_{0}^{t} \gamma(x,s) \dd s$.

Note that, $\Upsilon_{\lambda}$ is coercive and sequentially wekly lower semicontinuous. So, 
\[
\Upsilon_{\mu}(u_{\mu}) = \text{min}\left[\Upsilon_{\mu}(u); \ u \in W_0^{s_1,p}(\Omega)\right].
\]
is attaned by a function $u_{\mu}:= u_{k,\mu} \in W_0^{s_1,p}(\Omega)$.

Thus, $\Upsilon_{\mu}'(u_{\mu}) = 0 $, that is,
\begin{equation}\label{test2}
\bigg\langle A_{s_1,p}(u_{\mu}) + A_{s_2,q}(u_{\mu}),\varphi \bigg\rangle = \int_{\Omega}\gamma_{\mu}(x,u_{\mu}) \varphi \dd x, \ \ \text{for all} \ \ \varphi \in W_{0}^{s_1,p}(\Omega).
\end{equation}

Testing the equation \eqref{test2} with $\varphi = (u_\mu - u_{\lambda})^{+} \in W_{0}^{s_1,p}(\Omega)$, using the Proposition \ref{PC} and $0<\mu<\lambda$ we show that $u_{\mu} \leq u_{\lambda}$. In addition, testing the equation \eqref{test2} with the function $\varphi = (\tilde{u}_{\mu} - u_{\mu})^+ \in W_0^{s_1,p}(\Omega)$, using the Proposition \ref{PC} and the fact $\tilde{u}_{\mu}$ is unique solution of the problem \eqref{P4}, we show $\tilde{u}_{\mu} \leq u_{\mu}$.

So we have proved that,
\[
u_{\mu} \in [\tilde{u}_{\mu}, u_{\lambda}] \Rightarrow u_{\mu} \in S_{\mu} \subseteq \text{int}[(C_{s_1}^{0}(\overline{\Omega}))_+] \ \ \text{and so} \ \ \mu \in \mathcal{L}.
\]
\qed
\begin{propo}\label{prop5}
If hypotheses $(\textbf{H}_k)$ hold, $\lambda \in \mathcal{L}$, $u_{\lambda} \in S_{\lambda} \subseteq \text{int}\left[(C_{s_1}^{0}(\overline{\Omega}))_+\right]$ and $\mu < \lambda$, then $\mu \in \mathcal{L}$ and there exists $u_\mu \in S_\mu$ such that
\[
u_\lambda - u_\mu \in \text{int}\left[(C_{s_1}^{0}(\overline{\Omega}))_+\right].
\]
\end{propo}
\n {\bf Proof.} By Proposition \ref{propex} we know that $\mu \in \mathcal{L}$ and we can find $u_\mu:= u_{k,\mu} \in S_{\mu} \subseteq \text{int}\left[(C_{s_1}^{0}(\overline{\Omega})_+)\right]$ such that $u_\mu \leq u_\lambda$. Let $\rho = \Vert u_{\lambda} \Vert_{\infty}$ and $\widehat{E}_{k,\rho} > 0$ be as postulated by hypothesis $(\textbf{H}_k) \ (v)$. We have
\begin{align*}
(-\Delta_p)^{s_1} u_{\mu}(x) &+ (-\Delta_q)^{s_2} u_{\mu}(x) + \lambda\widehat{E}_{k,\rho} u_{\mu}(x)^{p-1} - \lambda u_{\mu}(x)^{-\eta} \\
&\leq \mu f_k(x,u_{\mu}(x)) + \lambda\widehat{E}_{k,\rho} u_{\mu}(x)^{p-1}\\
&= \lambda \left[f_k(x,u_{\mu}(x))+\widehat{E}_{k,\rho} u_{\mu}(x)^{p-1}\right] - \left(\lambda - \mu \right)f_k(x,u_{\mu}(x)) \\
&\leq \lambda \left[f_k(x,u_{\mu}(x)) + \widehat{E}_{k,\rho} u_{\mu}(x)^{p-1}\right]\\
&= (-\Delta_p)^{s_1} u_{\lambda}(x) + (-\Delta_q)^{s_2} u_{\lambda}(x)  + \lambda\widehat{E}_{k,\rho} u_{\lambda}(x)^{p-1} - \lambda u_{\lambda}(x)^{-\eta}.
\end{align*}

Note that, the function $g(t) = \lambda \widehat{E}_{k,\rho}t^{p-1} - \lambda t^{-\eta}$ is nondecreasing in $\mathbb{R}_{0}^{+}$, thus, by Proposition \ref{SPC} we have $u_\lambda - u_\mu \in \text{int}\left[(C_{s_1}^{0}(\overline{\Omega}))_+\right].$
\qed
\begin{propo}\label{prop4}
Assume that the hypotheses $(\textbf{H}_k)$ hold. Then $\lambda^{*} = \sup \mathcal{L} < +\infty,$ for each $k \in \mathbb{N}$.
\end{propo}
\n {\bf Proof.} By hyphoteses $H(i),(ii)$ and $(iii)$ we can find $\widehat{\lambda} > 0$ such that
\begin{equation}\label{eq4}
t^{p-1} \leq \widehat{\lambda} f_k(x,t) \ \ \text{for all} \ \ x \in \Omega, \ \ \text{all} \ \ t \geq 0.
\end{equation}

Let $\lambda>\lambda^{*}$ and suppose that $\lambda \in \mathcal{L}$. Then, there exists $u_\lambda:= u_{k,\lambda} \in S_{\lambda} \subseteq \text{int}[(C_{s_1}^{0}(\overline{\Omega}))_+]$, that is, $u_\lambda$ is a solution of the problem \eqref{P1}. Consider $\Omega_0 \subset\subset \Omega$ and $m_0 = \displaystyle\min_{\overline{\Omega}}u_{\lambda} > 0$. For $\delta \in (0,1)$ small we set $m_0^\delta = m_0 + \delta$. Let $\rho = \Vert u_\lambda \Vert_{\infty}$ and $\widehat{E}_{k,\rho} > 0$ be as postulated by $H(v)$. We have,
\begin{align*}
(-\Delta_p)^{s_1} m_0^{\delta} &+ (-\Delta_q)^{s_2} m_0^{\delta} + \lambda\widehat{E}_{k,\rho} (m_0^\delta)^{p-1} - \lambda (m_0^\delta)^{-\eta} \\
&\leq \lambda\widehat{E}_{k,\rho} (m_0^\delta)^{p-1} + \chi(\delta) \ \ (\text{with} \ \ \chi(\delta) \to 0^+ \ \ \text{as} \ \ \delta \to 0^+)\\
&= \left[\lambda\widehat{E}_{k,\rho} + 1\right] m_0^{p-1} + \chi(\delta) \\
&\leq \widehat{\lambda} f_k(x, m_0) + \lambda\widehat{E}_{k,\rho} (m_0^\delta)^{p-1} + \chi(\delta) \ \ (\text{see} \ \eqref{eq4})\\
&= \lambda \left[f_k(x,m_0)  + \widehat{E}_{k,\rho} (m_0^\delta)^{p-1}\right] - (\lambda - \widehat{\lambda}) f_k(x, m_0)+ \chi(\delta)\\
&\leq \lambda\left[f_k(x,u_\lambda(x)) + \widehat{E}_{k,\rho} u_{\lambda}^{p-1}\right] \ \ \text{for} \ \delta (0,1) \ \ \text{small enough}.\\
&= (-\Delta_p)^{s_1} u_{\lambda}(x) + (-\Delta_q)^{s_2} u_{\lambda}(x)  + \lambda\widehat{E}_{k,\rho} u_{\lambda}(x)^{p-1} - \lambda u_{\lambda}(x)^{-\eta}.
\end{align*}
where we have used the hyphoteses $H(iv), (v)$ and the fact $\chi(\delta) \to 0^{+}$ as $\delta \to 0^{+}$. By strong comparison principle we have 
\[
u_\lambda - m_0^{\delta} \in \text{int}[(C_{s_1}^{0}(\Omega))_+] \ \ \text{for} \ \ \delta \in (0,1) \ \ \text{small enougth}
\]
which contradicts with the definition of $m_0$. Consequently, it holds $0<\lambda^* \leq \widehat{\lambda} < \infty$. 
\qed
\begin{propo}\label{propo1}
If hyphoteses $(\textbf{H}_k)$ hold and $\lambda \in (0,\lambda^*)$, then problem $\eqref{P1}$ has least two positive solutions
\[
u_0, \hat{u} \in \text{int}[(C_{s_1}^{0}(\Omega))_+] \ \ \text{with} \ \ u_0 \leq \hat{u} \ \ \text{and} \ \ u_0 \neq \hat{u}.
\]
\end{propo}
\n {\bf Proof.} Let $0<\lambda < \vartheta < \lambda^*$. By Proposition \ref{prop4} $\lambda,\vartheta \in \mathcal{L}$. Thus, by Proposition \ref{prop5} we can find $u_0 \in S_{\lambda} \subseteq \text{int}[(C_{s_1}^{0}(\Omega))_+]$ and $u_{\vartheta} \in S_{\vartheta} \subseteq \text{int}[(C_{s_1}^{0}(\Omega))_+]$ such that
\[
u_{\vartheta} - u_0 \in S_{\lambda} \subseteq \text{int}[(C_{s_1}^{0}(\Omega))_+].
\]

From Proposition \ref{prop5}, we know that $u_\lambda \leq u_0$, hence $u_0^{-\eta} \in L^{1}(\Omega)$. Consider the Carathéodory function
\[
\widehat{\omega}_{\lambda}(x,t) = \left\{
\begin{array}{llll}
	\lambda[{u}_{0}^{-\eta} + f_k(x, {u}_{0})] & {\rm if} \ \ t < {u}_{0}(x),\\
	\lambda[t^{-\eta} + f_k(x, t)] & {\rm if} \ \  {u}_{0}(x) \leq t \leq {u}_{\vartheta}(x),\\
	\lambda[{u}_{\vartheta}^{-\eta} + f_k(x,{u}_{\vartheta})] & {\rm if} \ \ {u}_{\vartheta}(x)< t
\end{array}
\right.
\]
and define the $C^1$-functional $\widehat{\mu}_{\lambda}: W_{0}^{s_1,p}(\Omega) \to \mathbb{R}$ by
\[
\widehat{\mu}_{\lambda}(u) = \frac{1}{p} \left[ u \right]_{s_1,p}^{p} + \left[ u \right]_{s_2,p}^{q} - \int_{\Omega} \widehat{W}_{\lambda}(x,u) dx \ \ \text{for all} \ \ u \in W_{0}^{s_1,p}(\Omega).
\]
where $\widehat{W}_{\lambda}(t,x) = \displaystyle\int_{0}^{t} \widehat{\omega}_\lambda(x,s) \dd s$.

Consider also another Carathéodory function
\[
\omega_{\lambda}(x,t) = \left\{
\begin{array}{llll}
	\lambda[{u}_{0}^{-\eta}(x) + f_k(x, {u}_{0})] & {\rm if} \ \ t \leq {u}_{0}(x),\\
	\lambda[t^{-\eta} + f_k(x, t)] & {\rm if} \ \  {u}_{0}(x) < t
\end{array}
\right.
\]
and define the $C^1$-functional $\mu_{\lambda}: W_{0}^{s_1,p}(\Omega) \to \mathbb{R}$ by
\[
\mu_{\lambda}(u) = \frac{1}{p} \left[ u \right]_{s_1,p}^{p} + \left[ u \right]_{s_2,p}^{q} - \int_{\Omega} W_{\lambda}(x,u) dx \ \ \text{for all} \ \ u \in W_{0}^{s_1,p}(\Omega)
\]
where $W_{\lambda}(t,x) = \displaystyle\int_{0}^{t}\omega_\lambda (x,s) \dd s$.

It is clear that,
\begin{equation}\label{equ1}
\widehat{\mu}_{\lambda}(u)\bigg\vert_{[0,u_{\theta}]} = \mu_{\lambda}(u)\bigg\vert_{[0,u_{\theta}]} \ \ \text{and} \ \ \widehat{\mu}'_{\lambda}(u)\bigg\vert_{[0,u_{\theta}]} = \mu'_{\lambda}(u)\bigg\vert_{[0,u_{\theta}]}
\end{equation}
Let $K_\mu = \big\{u \in W_{0}^{s_1,p}(\Omega); \ \ \mu'(u) = 0\big\}$. Using the same arguments used in [\cite{Papa}, Proposition $8$] we can show that
\begin{align}
K_{\widehat{\mu}_{\lambda}} &\subseteq [u_0,u_\theta] \cap \text{int}[(C_{s_1}^{0}(\overline{\Omega}))_+]\label{equ2}\\
K_{\mu_{\lambda}} &\subseteq [u_0) \cap \text{int}[(C_{s_1}^{0}(\overline{\Omega}))_+]\label{equ3}
\end{align}

From \eqref{equ3}, we can assume that $K_{\mu_{\lambda}}$ is finite. Otherwise, we already have an infinity of positive smooth solutions of \eqref{P1} bigger than $u_0$ and so we are done. In adittion, we can assume that 
\begin{equation}\label{equ4}
K_{\mu_{\lambda}} \cap [u_0,u_{\theta}] = \{u_0\}.
\end{equation}

Moreover, it is clear that $\widehat{\mu}_{\lambda}$ is coercive and sequentially weakly lower semicontinuous. So there exists $\tilde{u}_0 \in W_{0}^{s_1,p}(\Omega)$ such that,
\[
\widehat{\mu}_{\tilde{u}_{0}} = \min\bigg[\widehat{\mu}_{\lambda}(u); \ \ u \in W_{0}^{s_1,p}(\Omega)\bigg]
\]
from \eqref{equ2} we have 
\[
\tilde{u}_{0} \in K_{\widehat{\mu}_{\lambda}} \subseteq [u_0,u_{\theta}] \cap \text{int}[(C_{s_1}^{0}(\overline{\Omega}))_+] 
\]
and so, from \eqref{equ1} and \eqref{equ4} results $\tilde{u}_{0} = u_0$. Therefore,
\[
u_0 \in \text{int}[(C_{s_1}^{0}(\overline{\Omega}))_+] \ \text{is a local} \ \ W_{0}^{s_1,p}(\Omega)-\text{minimizer of} \ \ \mu_{\lambda}.
\]

Consequently, there exists $\rho \in (0,1)$ such that,
\[
\mu_{\lambda}(u_0) < \inf\bigg[\mu_{\lambda}(u); \ \ \left[ u-u_0 \right]_{s_1,p} = \rho \bigg] = m_{\lambda}.
\]

Note that, if $u \in \text{int}[(C_{s_1}^{0}(\overline{\Omega}))_+]$, then on account of hypothesis $(\textbf{H}_k \ (ii))$ we have,
\[
\mu_{\lambda}(tu) \to -\infty \ \ \text{as} \ \ t\to \infty
\]
and moreover, classical arguments, which can be found in (\cite{Papa}, \cite{Sing}), along with conditions $(\textbf{H}_k)$ show that the function $\mu_{\lambda}$ satisfies the Cerami condition. By mountain pass theorem, there exists $\widehat{u} \in W_{0}^{s_1,p}(\Omega)$ such that,
\[
\widehat{u} \in K_{\mu_{\lambda}} \subseteq [u_0) \cap \text{int}[(C_{s_1}^{0}(\overline{\Omega}))_+]
\]
and $m_\lambda \leq \mu_{\lambda}(\widehat{u})$. So, we have $\widehat{u} \in S_{\lambda}, \ u_0 \leq \widehat{u}$ and $\widehat{u} \neq u_0$.
\qed
\begin{propo}
If hypotheses $(\textbf{H}_k)$ hold, then $\lambda^* \in \mathcal{L}$.
\end{propo}
\n {\bf Proof.} Let $\{\lambda_n\} \subset (0,\lambda^*)$ be such that  $\lambda_n \to \lambda^*$. We have $\{\lambda_n\}_{n\geq 1} \subseteq\mathcal{L}$ and of the proof of Proposition \ref{propo1} we find $u_n \in S_{\lambda_n} \subseteq \text{int}[(C_{s_1}^{0}(\overline{\Omega}))_+]$ such that,
\begin{align*}
\mu_{\lambda_n}(u_n) &= \frac{1}{p} \left[ u_n \right]_{s_1,p}^{p} + \left[ u_n \right]_{s_2,q}^{q} - \lambda_n\int_{\Omega} [u_n^{1-\eta} + f_k(x,u_n).u_n] \dd x \\
&= \frac{1}{p} \left[ u_n \right]_{s_1,p}^{p} + \frac{1}{q} \left[ u_n \right]_{s_2,q}^{q} - \left[ u_n \right]_{s_1,p}^{p} - \left[ u_n \right]_{s_2,p}^{p} \ \ ( \text{Since} \ \ u_n \in S_{\lambda_n})\\
&= \left(\frac{1}{p} - 1\right) \left[u_n \right]_{s_1,p}^{p} + \left(\frac{1}{q} - 1\right) \left[ u_n \right]_{s_2,q}^{q} < 0 \ \ \text{for all} \ \ n \in \mathbb{N}.
\end{align*}

Moreover, we have
\begin{equation}\label{test3}
\bigg\langle A_{s_1,p}(u_n) + A_{s_2,q}(u_n),\varphi \bigg\rangle = \int_{\Omega}[\lambda_n u_n^{-\eta} + f_k(x,u_n)]\varphi\dd x, \ \ \text{for all} \ \ \varphi \in W_{0}^{s_1,p}(\Omega).
\end{equation}

Arguing as in the proof of Proposition $13$ in \cite{Sing}, we obtain that at least for a subsequence, 
\[
u_n \to u_* \ \ \text{in} \ \ W_{0}^{s_1,p}(\Omega) \ \ \text{as} \ n \to \infty.
\]
By Proposition \ref{prop5}, $\tilde{u}_{\lambda_1} \leq u_n$ for all $n \in \mathbb{N}$. Therefore, we see $u_* \neq 0$ and $u_*^{-\eta}\varphi \leq \tilde{u}_{\lambda_1}^{-\eta}\varphi \in L^1(\Omega)$ for all $\varphi \in W_{0}^{s_1,p}(\Omega)$. In \eqref{test3}, we pass to the limit as $n \to \infty$ and we obtain
\begin{equation}\label{test3}
\bigg\langle A_{s_1,p}(u_*) + A_{s_2,q}(u_*),\varphi \bigg\rangle = \int_{\Omega}[\lambda^* u_*^{-\eta} + f_k(x,u_*)]\varphi\dd x, \ \ \text{for all} \ \ \varphi \in W_{0}^{s_1,p}(\Omega).
\end{equation}
that is,
\[
u_* \in S_{\lambda^*} \subseteq \text{int}[(C_{s_1}^{0}(\overline{\Omega}))_+] \ \ \text{and so} \ \ \lambda^* \in \mathcal{L}.
\]
\qed
So, summarizing the situation for problem \eqref{P1}, we can state the following bifurcation-type theorem.
\begin{teor}\label{teo1}
If hypotheses $(\textbf{H}_k)$ hold, then we can find $\lambda^*>0$ such that
\begin{enumerate}
\item for every $\lambda \in (0,\lambda^*)$ problem \eqref{P1} has at least two nontrivial positive solutions 
\[u_0, \hat{u} \in \text{int}[(C_{s_1}^{0}(\Omega))_+] \ \ \text{with} \ \ u_0 \leq \hat{u} \ \ \text{and} \ \ u_0 \neq \hat{u}.
\]
\item for $\lambda = \lambda^*$ problem \eqref{P1} has one nontrivial positive solution 
\[
u_* \in \text{int}[(C_{s_1}^{0}(\overline{\Omega}))_+] \ \ \text{and so} \ \ \lambda^* \in \mathcal{L}.
\]
\item for $\lambda > \lambda^*$ problem \eqref{P1} has no nontrivial positive solution.
\end{enumerate}

\end{teor}

\section{Existence of positive solution for \ref{P0}}

We denote by $u:= u_{k,\lambda}$ the solution of the problem (\ref{P1}) given by Theorem \ref{teo1}. Thus, we obtain

\begin{propo}\label{Linfty}
Let $u:= u_{k,\lambda} \in W_{0}^{s_1,p}(\Omega)$ be a positive weak solution to the problem in (\ref{P1}), then $u \in L^{\infty}(\bar{\Omega})$. Moreover, there exists $k >1$ sufficiently large such that,
\[
\Vert u \Vert_\infty \leq M_k.
\]
\end{propo}

\n {\bf Proof.} The arguments of the proof is taken from the celebrated article of \cite{IterMos} with appropriate modifications. We will proceed with the smooth, convex and Lipschitz function $g_\epsilon(t)=\left(\epsilon^2+t^2\right)^{\frac{1}{2}}$ for every $\epsilon>0$. Moreover, $g_\epsilon(t) \rightarrow|t|$ as $t \rightarrow 0$ and $\left|g_\epsilon^{\prime}(t)\right| \leq 1$. Let $0<\psi \in C_c^{\infty}(\Omega)$ and choose $\varphi=\psi\left|g_\epsilon^{\prime}(u)\right|^{p-2} g_\epsilon^{\prime}(u)$ as the test function.

By Lemma $5.3$ of \cite{IterMos} for all $\psi \in C_c^{\infty}(\Omega) \cap \mathbb{R}^{+}$, we obtain
$$
\begin{aligned}
\langle A_{s_1,p}(g_\epsilon(u)), \psi \rangle + \langle A_{s_2,q}(g_\epsilon(u)), \psi \rangle \leq \lambda \int_{\Omega}\left(\frac{1}{|u|^{\eta} }+ | f_k(x,u)| \right)\left|g'_\epsilon(u)\right|^{p-1} \psi \mathrm{d} x
\end{aligned}
$$
By Fatou's Lemma as $\varepsilon \to 0$ we have,
\begin{equation} \label{solu}
\begin{aligned}
\langle A_{s_1,p}(u), \psi \rangle + \langle A_{s_2,q}(u), \psi \rangle \leq \lambda \int_{\Omega}\left(\frac{1}{|u|^{\eta} }+ | f_k(x,u)| \right)\psi \mathrm{d} x
\end{aligned}
\end{equation}

Define $u_n = \min\{(u-M_k^{\gamma})^{+}, n\}$ for each $n \in \mathbb{N}$ and $\gamma >0$. Let $\beta > 1$ , $\delta > 0$ and consider $\psi_\delta = (u_n + \delta)^{\beta} - \delta^{\beta}$. Thus, $\psi_\delta = 0$ in $\{u\leq M_k^{\gamma}\}$ and using $\psi_\delta$ in \eqref{solu} we obtain

$$
\begin{aligned}
\langle A_{s_1,p}(u), \psi_\delta \rangle + \langle A_{s_2,q}(u), \psi_\delta \rangle \leq \lambda \int_{\Omega}\left(\frac{1}{|u|^{\eta} }+ | f_k(x,u)| \right)((u_n + \delta)^\beta - \delta^{\beta}) \mathrm{d} x
\end{aligned}
$$

By Lema $5.4$ in \cite{IterMos} to follow the estimates,

\begin{equation*}
\begin{aligned}
\langle A_{s_1,p}(u), \psi_\delta \rangle + \langle A_{s_2,q}(u), \psi_\delta \rangle &\geq \beta \left(\displaystyle\frac{p}{\beta + p -1} \right)^{p} \left[(u_n + \delta)^{\frac{\beta + p -1}{p}}\right]_{s_1,p}^{p} + \beta \left(\displaystyle\frac{q}{\beta + q -1} \right)^{q} \left[(u_n + \delta)^{\frac{\beta + q -1}{q}}\right]_{s_2,q}^{q}\\
&\geq \beta \left(\displaystyle\frac{p}{\beta + p -1} \right)^{p} \left[(u_n + \delta)^{\frac{\beta + p -1}{p}}\right]_{s_1,p}^{p}
\end{aligned}
\end{equation*}
consequently,

\begin{equation*}
\begin{aligned}
\beta \left(\displaystyle\frac{p}{\beta + p -1} \right)^{p} \left[(u_n + \delta)^{\frac{\beta + p -1}{p}}\right]_{s_1,p}^{p} \leq \lambda \int_{\Omega}\left(\frac{1}{|u|^{\eta} }+ | f_k(x,u)| \right)((u_n + \delta)^\beta - \delta^{\beta}) \mathrm{d} x
\end{aligned}
\end{equation*}
and thus,
\begin{equation} \label{estA}
\begin{aligned}
\left[(u_n + \delta)^{\frac{\beta + p -1}{p}}\right]_{s_1,p}^{p} \leq \lambda \displaystyle\frac{1}{\beta} \left(\displaystyle\frac{\beta + p -1}{p} \right)^{p} \int_{\Omega}\left(\frac{1}{|u|^{\eta} }+ | f_k(x,u)| \right)((u_n + \delta)^\beta - \delta^{\beta}) \mathrm{d} x
	\end{aligned}
\end{equation}

Using the estimates \eqref{est.1}, for $M_k > 1$ we have,
$$
\begin{aligned}
\int_{\Omega}\left(\frac{1}{|u|^{\eta} }+ | f_k(x,u)| \right)&((u_n + \delta)^\beta - \delta^{\beta}) \mathrm{d} x \leq  \int_{\Omega}\left(\frac{1}{|u|^{\eta} } + C.M_k^{2\theta} \vert u \vert^{r-1}\right)((u_n + \delta)^\beta - \delta^{\beta}) \mathrm{d} x \\
&= \int_{\{u \geq M_k^\gamma\}}\left(\frac{1}{|u|^{\eta} } + C.M_k^{2\theta} \vert u \vert^{r-1}\right)((u_n + \delta)^\beta - \delta^{\beta}) \mathrm{d} x \\
&= \int_{\{u \geq M_k^\gamma \}} ((u_n + \delta)^\beta - \delta^{\beta}) \mathrm{d} x + \int_{\{u \geq M_n \}} C.M_k^{2\theta} \vert u \vert^{r-1}((u_n + \delta)^\beta - \delta^{\beta}) \mathrm{d} x \\
&\leq \int_{\{u \geq M_k^\gamma\}} M_k^{2\theta} ((u_n + \delta)^\beta - \delta^{\beta}) \mathrm{d} x   + \int_{\{u \geq M_k^\gamma \}} C.M_k^{2\theta} \vert u \vert^{r-1}((u_n + \delta)^\beta - \delta^{\beta}) \mathrm{d} x \\
& \leq \int_{\Omega} M_k^{2\theta} ((u_n + \delta)^\beta - \delta^{\beta}) \mathrm{d} x   + \int_{\Omega} C.M_k^{2\theta} \vert u \vert^{r-1}((u_k + \delta)^\beta - \delta^{\beta}) \mathrm{d} x\\
& \leq C.M_k^{2\theta} \left(\vert \Omega \vert^{\frac{\sigma-1}{\sigma}}  + \Vert u \Vert^{r-1}_{L^{p_{s_1}^*}(\Omega)} \right) \Vert (u_n + \delta)^{\beta} \Vert_{L^{\sigma}(\Omega)}
\end{aligned}
$$
where $C$ is a constant independent of $k$ and $\sigma = \displaystyle\frac{p_{s_1}^{*}}{p_{s_1}^{*}-r+1}$. Morevover, observe that, the function $u:= u_k$ satisfies $u \leq \overline{u}$ where $\overline{u}$ is a supersolution of the problema \eqref{P3} does not depend on $k$, we have $\Vert u \Vert^{r-1}_{L^{p_{s_1}^*}(\Omega)} \leq C_0 \Vert \overline{u} \Vert^{r-1}_{\infty}$ independent of $k$. Thus,

\begin{equation}\label{estB}
\begin{aligned}
\int_{\Omega}\left(\frac{1}{|u|^{\eta} }+ | f_k(x,u)| \right)&((u_n + \delta)^\beta - \delta^{\beta}) \mathrm{d} x \leq
K M_k^{2\theta} \left(\vert \Omega \vert^{\frac{\sigma-1}{\sigma}}  + \Vert \overline{u} \Vert^{r}_{\infty} \right) \Vert (u_n + \delta)^{\beta} \Vert_{L^{\sigma}(\Omega)}\\
& = K_0 M_k^{2 \theta} \Vert (u_n + \delta)^{\beta} \Vert_{L^{\sigma}(\Omega)}
\end{aligned}
\end{equation}
with $K_0$ independent of $k$.

By Sobolev inequality, triangle inequality and $(u_n + \delta)^{\beta + p -1} \geq \delta^{p-1} (u_n + \delta)^{\beta}$
\begin{equation}\label{estC}
\begin{aligned}
\left[(u_n + \delta)^{\frac{\beta + p -1}{p}}\right]_{s_1,p}^{p} &\geq S \Vert (u_n+ \delta)^\beta - \delta^{\beta} \Vert_{L^{p_{s_1}^{*}}(\Omega)}^p \\
&\geq \left(\displaystyle\frac{\delta}{2}\right)^{p-1} \left[ \int_{\Omega} \vert (u_n+\delta)^{\frac{p_{s_1}^{*}\beta}{p}} \dd x \right]^{\frac{p}{p_{s_1}^{*}}} - \delta^{\beta + p-1} \vert \Omega \vert^{\frac{p}{p_{s_1}^{*}}}\\
&\geq \left(\displaystyle\frac{\delta}{2}\right)^{p-1} \Vert (u_n + \delta)^{\frac{\beta}{p}} \Vert_{L^{p_{s_1}^{*}}(\Omega)}^{p} - M_k^{2\theta} \delta^{\beta+p-1} \vert \Omega \vert^{\frac{p}{p_{s_1}^{*}}},
\end{aligned}
\end{equation}
in the estimate above we using that $M_k > 1$.

Using thes estimates \eqref{estC} and \eqref{estB} in \eqref{estA}  we obtain,

$$
\begin{aligned}
\left\|\left(u_n+\delta\right)^{\frac{\beta}{p}}\right\|_{L^{p_{s_1}^{*}}(\Omega)}^p &\leq \left(\frac{2}{\delta}\right)^{p-1}\left[\left(\frac{ (\beta+p-1)^p}{\beta p^p}\right) K_0 M_k^{2\theta} \Vert (u_n + \delta)^{\beta} \Vert_{L^{\sigma}(\Omega)} +\delta^{\beta+p-1}|\Omega|^ \frac{p}{p_s^*}\right] \\
&= \left(\frac{2}{\delta}\right)^{p-1} \left(\frac{ (\beta+p-1)^p}{\beta p^p}\right) K_0 M_k^{2\theta} \Vert (u_n + \delta)^{\beta} \Vert_{L^{\sigma}(\Omega)} +\delta^{\beta}|\Omega|^ \frac{p}{p_s^*}\\
&\leq \left(\frac{2}{\delta}\right)^{p-1} \left(\frac{ (\beta+p-1)^p}{\beta p^p}\right) K_0 M_k^{2\theta} \Vert (u_n + \delta)^{\beta} \Vert_{L^{\sigma}(\Omega)} + |\Omega|^ {\frac{p}{p_s^*} - 1} \int_{\Omega} (u_n + \delta)^{\beta}  \dd x
\end{aligned}
$$

By Holder's inequality, we have 
\[
\delta^{\beta} = \vert \Omega \vert^{-1} \int_{\Omega} \delta^{\beta} \dd x \leq \vert \Omega \vert^{-1} \int_{\Omega} (u_n + \delta)^{\beta} \dd x \leq \vert \Omega \vert^{-\frac{1}{\sigma}} \Vert (u_n + \delta)^{\beta} \Vert_{L^{\sigma}(\Omega)}.
\]

Consequently, 
$$
\begin{aligned}
\left\|\left(u_n+\delta\right)^{\frac{\beta}{p}}\right\|_{L^{p_{s_1}^{*}}(\Omega)}^p &\leq \left(\frac{2}{\delta}\right)^{p-1} \left(\frac{ (\beta+p-1)^p}{\beta p^p}\right) K_0 M_k^{2\theta} \Vert (u_n + \delta)^{\beta} \Vert_{L^{\sigma}(\Omega)} + |\Omega|^ {\frac{p}{p_s^*} - \frac{1}{\sigma}} \Vert (u_n + \delta)^{\beta} \Vert_{L^{\sigma}(\Omega)}.\\
\end{aligned}
$$
Since, $\displaystyle\frac{1}{\beta}\left(\frac{\beta+p-1}{p}\right)^p \geq 1$ we can deduce that
$$
\begin{aligned}
\left\|\left(u_n+\delta\right)^{\frac{\beta}{p}}\right\|_{L^{p_{s_1}^{*}}(\Omega)}^p &\leq \frac{1}{\beta}\left(\frac{\beta+p-1}{p}\right)^p M_k^{2\theta} \left\|\left(u_n+\delta\right)^\beta\right\|_q\left(\frac{K_0}{\delta^{p-1}}+|\Omega|^{\frac{p}{p_s^*}-\frac{1}{\sigma}}\right)\\
\end{aligned}
$$

Now choose, $\delta > 0$ such that $\delta^{p-1} = K_0 \vert \Omega \vert^{\frac{1}{\sigma} - \frac{p}{p_{s_1}^{*}}}$ and $\beta >1$ such that, $\left(\displaystyle\frac{\beta+p-1}{p}\right)^p \geq \beta^{p}$. Thus, 
$$
\begin{aligned}
\left\|\left(u_n+\delta\right)^{\frac{\beta}{p}}\right\|_{L^{p_{s_1}^{*}}(\Omega)}^p \leq C M_k^{2\theta} \beta^{p-1} \left\|\left(u_n+\delta\right)^\beta\right\|_{L^{\sigma}(\Omega)}\\
\end{aligned}
$$

For $\tau = \sigma \beta$ and $\alpha =\displaystyle\frac{p_{s_1}^{*}}{\sigma p}$ we obtain,
$$
\begin{aligned}
\left\|u_n+\delta\right\|_{L^{\alpha \tau}(\Omega)}^\beta \leq C M_k^{2\theta} \beta^{p-1} \left\|u_n+\delta \right\|_{L^{\tau}(\Omega)}^\beta\\
\end{aligned}
$$
and therefore,
$$
\begin{aligned}
\left\|u_n+\delta\right\|_{L^{\alpha \tau}(\Omega)} \leq \left( C M_k^{2\theta}\right)^{\frac{\sigma}{\tau}} \left(\displaystyle\frac{\tau}{\sigma}\right)^{(p-1)\frac{\sigma}{\tau}} \left\|u_n+\delta \right\|_{L^{\tau}(\Omega)}.
\end{aligned}
$$

Taking, $\tau_0 = \sigma$, $\tau_{m+1} = \alpha \tau_m = \alpha^{m+1} \sigma$, then after performing $m$ iterations we obtain the inequality

$$
\begin{aligned}
\left\|u_n+\delta\right\|_{L^{\tau_{m+1}}(\Omega)} &\leq \left( C M_k^{2\theta}\right)^{\displaystyle\sum_{i=0}^{m}\frac{\sigma}{\tau_i}} \left(\displaystyle\prod_{i=1}^{m}\left(\displaystyle\frac{\tau_i}{\sigma}\right)^{\frac{\sigma}{\tau_i}}\right)^{(p-1)} \left\|u_n+\delta \right\|_{L^{\tau}(\Omega)}\\
&= \left( C M_k^{2\theta}\right)^{\displaystyle\sum_{i=1}^{m}\frac{1}{\alpha^{i}}} \left(\displaystyle\prod_{i=1}^{m}\alpha^{\frac{i}{\alpha^i}}\right)^{(p-1)} \left\|u_n+\delta \right\|_{L^{\tau}(\Omega)}
\end{aligned}
$$

Therefore, on passing the limit as $m \to \infty$, we get
\begin{equation}\label{emb.result}
	\begin{aligned}
\Vert u_n \Vert_{L^{\infty}(\Omega)} \leq \left\|u_n+\delta\right\|_{L^{\infty}(\Omega)} \leq  C^{\frac{\alpha}{\alpha-1}} M_k^{\frac{2\theta \alpha}{\alpha-1}} \alpha^{\frac{(p-1)\alpha}{(\alpha-1)^2}} \left\|u_n+\delta \right\|_{L^{\sigma}(\Omega)} \leq C_1 M_k^{\frac{2\theta \alpha}{\alpha-1}}.
\end{aligned}
\end{equation}

In the last inequality we use the fact, $u \leq \overline{u}$, where $\overline{u} \in L^{\infty}(\Omega)$ is a supersolution of the problem \eqref{P3} and thus, $u_n = \min\{(u - M_k^{\gamma})^{+}, n\} \leq (u-M_k^{\gamma})^+ \leq u^{+} \leq \overline{u}$, for each $n \in \mathbb{N}$ and $k$ large enough (such that $\Vert \overline{u} \Vert \leq M_k^{\gamma})$.

Therefore, as $n \to \infty$ we obtain
\[
\Vert (u - M_k^{\gamma})^+ \Vert_{\infty} \leq M_k
\]
for $M_k$ sufficiently large and $\displaystyle\frac{2\theta \alpha}{\alpha-1} < 1$. Consequently, since $M_k \to \infty$ as $k \to \infty$ we have, for $\gamma < 1$, there exists $k>1$ large enough such that, 
\[
\Vert u \Vert_{\infty} \leq M_k.
\]

Also, by \eqref{emb.result}, the embedding $W^{s_1,p}_0(\Omega) \hookrightarrow L^{\sigma}(\Omega)$ and since $u_n = \min\{(u - M_k^{\gamma})^{+}, n\} \leq (u-M_k^{\gamma})^+ \leq u^{+}\leq |u|$ we can establish 
\[
\|u_n\|_{L^{\infty}\left(\Omega\right)} \leq CM_k^{\frac{2\theta \alpha}{\alpha-1}}[u]_{s_1,p}.
\]
Therefore, as $n \to \infty$ we obtain
\[
\|u\|_{L^{\infty}\left(\Omega\right)} \leq CM_k^{\frac{2\theta \alpha}{\alpha-1}}[u]_{s_1,p},
\]
for $k>1$ large enough fixed.
\qed

\begin{teor}\label{teo2}
If hypotheses $(\textbf{H})$ hold, then we can find $\lambda^*=\lambda^*(k)>0$ ($k$ as in Proposition \ref{Linfty}) such that
\begin{enumerate}
\item for every $\lambda \in (0,\lambda^*)$ problem \eqref{P0} has at least two nontrivial positive solutions 
\[u_0, \hat{u} \in \text{int}[(C_{s_1}^{0}(\Omega))_+] \ \ \text{with} \ \ u_0 \leq \hat{u} \ \ \text{and} \ \ u_0 \neq \hat{u}.
\]
\item for $\lambda = \lambda^*$ problem \eqref{P0} has one nontrivial positive solution 
\[
u_* \in \text{int}[(C_{s_1}^{0}(\overline{\Omega}))_+] \ \ \text{and so} \ \ \lambda^* \in \mathcal{L}.
\]
\item for $\lambda > \lambda^*$ problem \eqref{P0} has no nontrivial positive solution.
\end{enumerate}
\end{teor}

\n {\bf Proof.} By Theorem \ref{teo1}, for each $\lambda \in (0,\lambda^{*}]$ and $k \in \mathbb{N}$ there exists $u_{k,\lambda}$ such that,
\begin{equation} \label{P5}\tag{$P_{k,\lambda}$}
\left\{
\begin{array}{llll}
(-\Delta_p)^{s_1}u + (-\Delta_q)^{s_2}u = \lambda \left[u(x)^{-\eta} + f_k(x,u)\right] & {\rm in} \ \ \Omega,\\
u= 0 & {\rm in} \ \R^N\setminus\Omega,\\
u> 0 & {\rm in} \ \Omega.
\end{array}
\right.
\end{equation}
Moreover, $1,2$ and $3$ holds to the problem \eqref{P1}, by Theorem \ref{teo1}.

Using the Proposition \ref{Linfty}, we have $\Vert u_{k,\lambda} \Vert_{\infty} < M_k$ for some $k > 1$ large enough. Thus, $u_\lambda:= u_{k,\lambda}(x) \leq M_k$ and therefore $f_k(x,u_\lambda) = f(x,u_\lambda)$, in other words $u_\lambda$ satisfies the problem \eqref{P0}.
\qed

\noindent \textbf{Disclosure statement:}
On behalf of all authors, the corresponding author states that there is no conflict of interest.

\noindent \textbf{Funding}
The first author was partially supported by FAPEMIG/ APQ-02375-21, FAPEMIG/RED00133-21, and CNPq Processes 101896/2022-0, 305447/2022-0.

All the authors take part in the project projeto Special Visiting Researcher - FAPEMIG CEX APQ 04528/22.

\bibliographystyle{amsplain} 

\end{document}